\documentclass[12pt,thmsa]{article}
\usepackage{amssymb}
\usepackage{amsfonts}
\usepackage{amsmath}

\numberwithin{equation}{section}
\newtheorem{theorem}{Theorem}[section]

\newtheorem{lemma}[theorem]{Lemma}

\newtheorem{proposition}[theorem]{Proposition}

\begin{document}

\title{Multiplicity and concentration of nontrivial solutions for the
generalized extensible beam equations}
\date{}
\author{Juntao Sun$^{a},$\thanks{%
E-mail address: sunjuntao2008@163.com(J. Sun)}\quad Tsung-fang Wu$^b,$%
\thanks{%
E-mail address: tfwu@nuk.edu.tw (T.-F. Wu)} \\
{\footnotesize $^a$\emph{School of Science, Shandong University of
Technology, Zibo 255049, P.R. China }}\\
{\footnotesize $^b$\emph{Department of Applied Mathematics, National
University of Kaohsiung, Kaohsiung 811, Taiwan }}}
\maketitle

\begin{abstract}
In this paper, we study a class of generalized extensible beam equations
with a superlinear nonlinearity
\begin{equation*}
\left\{
\begin{array}{ll}
\Delta ^{2}u-M\left( \Vert \nabla u\Vert _{L^{2}}^{2}\right) \Delta
u+\lambda V(x) u=f( x,u) & \text{ in }\mathbb{R}^{N}, \\
u\in H^{2}(\mathbb{R}^{N}), &
\end{array}%
\right.
\end{equation*}%
where $N\geq 3$, $M(t) =at^{\delta }+b$ with $a,\delta >0$ and $b\in \mathbb{%
R}$, $\lambda >0$ is a parameter, $V\in C(\mathbb{R}^{N},\mathbb{R})$ and $%
f\in C(\mathbb{R}^{N}\times \mathbb{R},\mathbb{R}).$ Unlike most other
papers on this problem, we allow the constant $b$ to be nonpositive, which
has the physical significance. Under some suitable assumptions on $V(x)$ and
$f(x,u)$, when $a$ is small and $\lambda$ is large enough, we prove the
existence of two nontrivial solutions $u_{a,\lambda }^{(1)}$ and $%
u_{a,\lambda }^{(2)}$, one of which will blow up as the nonlocal term
vanishes. Moreover, $u_{a,\lambda }^{(1)}\rightarrow u_{\infty}^{(1)}$ and $%
u_{a,\lambda }^{(2)}\rightarrow u_{\infty}^{(2)}$ strongly in $H^{2}(\mathbb{%
R}^{N})$ as $\lambda\rightarrow\infty$, where $u_{\infty}^{(1)}\neq
u_{\infty}^{(2)}\in H_{0}^{2}(\Omega )$ are two nontrivial solutions of
Dirichlet BVPs on the bounded domain $\Omega$. It is worth noting that the
regularity of weak solutions $u_{\infty}^{(i)}(i=1,2)$ here is explored.
Finally, the nonexistence of nontrivial solutions is also obtained for $a$
large enough.
\end{abstract}


\section{Introduction}

Consider the nonlinear generalized extensible beam equations in the form:%
\begin{equation}
\left\{
\begin{array}{ll}
\Delta ^{2}u-M\left( \Vert \nabla u\Vert _{L^{2}}^{2}\right) \Delta
u+\lambda V\left( x\right) u=f\left( x,u\right) & \text{ in }\mathbb{R}^{N},
\\
u\in H^{2}(\mathbb{R}^{N}), &
\end{array}%
\right.  \tag{$E$}
\end{equation}%
where $N\geq 3,\Delta ^{2}u=\Delta (\Delta u),M\left( t\right) =at^{\delta
}+b$ with $a,\delta >0$ and $b\in \mathbb{R}$, $\lambda >0$ is a parameter,
and $f\in C(\mathbb{R}^{N}\times \mathbb{R},\mathbb{R}).$ We assume that the
potential $V(x)$ satisfies the following assumptions:

\begin{itemize}
\item[$(V1)$] $V\in C(\mathbb{R}^{N})$ and $V(x)\geq 0$ for all $x\in
\mathbb{R}^{N};$

\item[$(V2)$] there exists $c_{0}>0$ such that the set $\{V<c_{0}\}:=\{x\in
\mathbb{R}^{N}\ |\ V(x)<c_{0}\} $ has finite positive Lebesgue measure for $%
N\geq 4$ and
\begin{equation*}
\left\vert \{V<c_{0}\}\right\vert <S_{\infty }^{-2}\left( 1+\frac{A_{0}^{2}}{%
2}\right) ^{-1}\text{ for }N=3,
\end{equation*}%
where $\left\vert \cdot \right\vert $ is the Lebesgue measure, $S_{\infty }$
is the best Sobolev constant for the imbedding of $H^{2}(\mathbb{R}^{N})$ in
$L^{\infty }(\mathbb{R}^{N})$ for $N=3$, and $A_{0}$ is defined in (\ref{1-6}%
) below;

\item[$(V3)$] $\Omega =\text{int}\{x\in \mathbb{R}^{N}:V(x)=0\}$ is nonempty
and has smooth boundary with $\overline{\Omega}=\{x\in \mathbb{R}%
^{N}:V(x)=0\}.$
\end{itemize}

The hypotheses $(V1)-(V3),$ suggested by Bartsch et. al. \cite{BPW}, imply
that $\lambda V(x)$ represents a potential well whose depth is controlled by
$\lambda $. If $\lambda $ is sufficiently large, then $\lambda V(x)$ is
known as the steep potential well. About its applications, we refer the
reader to \cite{JZ,SW,SW1,SW2,SCW,SWW,YT1,ZLZ} and references therein.

Eq. $(E)$ arises in an interesting physical context. In 1950, Woinowsky and
Krieger \cite{W} introduced the following extensible beam equation:

\begin{equation}
\rho u_{tt}+EIu_{xxxx}-\left( \frac{Eh}{2I}\int\nolimits_{0}^{L}\left\vert
u_{x}\right\vert ^{2}dx+P_{0}\right) u_{xx}=0,  \label{1-4}
\end{equation}%
where $L$ is the length of the beam in the rest position, $E$ is the Young
modulus of the material, $I$ is the cross-sectional moment of inertia, $\rho$
is the mass density, $P_{0}$ is the tension in the rest position and $h$ is
the cross-sectional area. This model is used to describe the transverse
deflection $u(x,t)$ of an extensible beam of natural length $L$ whose ends
are held a fixed distance apart. Such problems are often referred to as
being nonlocal because of the presence of the term $\left(
\int\nolimits_{0}^{L}\left\vert u_{x}\right\vert ^{2}dx\right) u_{xx}$,
which indicates the change in the tension of the beam due to its
extensibility. The qualitative and stable analysis of solutions for Eq. (\ref%
{1-4}) can be traced back to the 1970s, for instance in the papers by Ball
\cite{B2}, Dickey \cite{D} and Medeiros \cite{M1}.

As a simplification of the von Karman plate equation, Berger \cite{B3}
proposed the plate model describing large deflection of plate as follows%
\begin{equation}
u_{tt}+\Delta ^{2}u-\left( \int_{\Omega}\left\vert \nabla u\right\vert
^{2}dx+Q_{0}\right) \Delta u=f\left( u,u_{t},x\right) ,  \label{1-5}
\end{equation}%
where $\Omega\subset \mathbb{R}^{N}(N=1,2)$ is a bounded domain with a
sufficiently smooth boundary, the parameter $Q_{0}$ is in-plane forces
applied to the plate ($Q_{0}>0$ represents outward pulling forces and $%
Q_{0}<0$ means inward extrusion forces) and the function $f$ represents
transverse loads which may depend on the displacement $u$ and the velocity $%
u_{t}$. Apparently, when $N=1$ and $f\equiv 0$ in Eq. (\ref{1-5}), the
corresponding equation becomes the extensible beam equation (\ref{1-4}).
Owing to its importance, the various properties of solutions for Eq. (\ref%
{1-5}) have been treated by many researchers; see for example, \cite%
{CCS,MN,P,Y,Z}. More precisely, Patcheu \cite{P} investigated the existence
and decay property of global solutions to the Cauchy problem of Eq. (\ref%
{1-5}) with $f\left( u,u_{t},x\right) \equiv f\left( u_{t}\right) $ in the
abstract form. Yang \cite{Y} studied the global existence, stability and the
longtime dynamics of solutions to the initial boundary value problem (IBVP)
of an extensible beam equation with nonlinear damping and source terms in
any space dimensions, i.e. Eq. (\ref{1-5}) with $f\left( u,u_{t},x\right)
=g(u_{t})+h(u)+k(x).$

In the last two decades, the stationary form of Eq. (\ref{1-5}), of the form
similar to Eq. $(E)$, has begun to attract attention, specially on the
existence and multiplicity of nontrivial solutions, but the relevant results
are rare. We refer the reader to \cite{CF,FN,LZ,M0,WAA,WAA1,XC,XC1} and
references therein. To be precise, Ma \cite{M0} studied the existence of
nontrivial solutions for a class of extensible beam equations with nonlinear
boundary conditions in dimension one. Wang et al. \cite{WAA} concentrated on
the following Navier BVPs:
\begin{equation}
\left\{
\begin{array}{ll}
\Delta ^{2}u+\lambda \left( a\int_{\Omega }\left\vert \nabla u\right\vert
^{2}dx+b\right) \Delta u=f\left( x,u\right) & x\in \Omega , \\
u=\Delta u=0 & x\in \partial \Omega ,%
\end{array}%
\right.  \label{1-1}
\end{equation}%
where $\Omega \subset \mathbb{R}^{N}$ is a smooth bounded domain and $%
\lambda ,a,b>0$. Applying mountain pass techniques and the truncation
method, they obtained the existence of nontrivial solutions for Eq. (\ref%
{1-1}) for $\lambda $ small enough when $f(x,u)$ satisfies some superlinear
assumptions. Cabada and Figueiredo \cite{CF} considered a class of
generalized extensible beam equations with critical growth in $\mathbb{R}%
^{N} $ as follows
\begin{equation}
\left\{
\begin{array}{ll}
\Delta ^{2}u-M\left( \Vert \nabla u\Vert _{L^{2}}^{2}\right) \Delta
u+u=\lambda f\left( u\right) +|u|^{2^{\ast \ast }-2}u & \text{ in }\mathbb{R}%
^{N}, \\
u\in H^{2}(\mathbb{R}^{N}), &
\end{array}%
\right.  \label{1-7}
\end{equation}%
where $M:\mathbb{R}^{+}\rightarrow \mathbb{R}^{+}$ are continuous increasing
functions, $f\in C(\mathbb{R},\mathbb{R})$, $2^{\ast \ast }=\frac{2N}{N-4}$
with $N\geq 5$ and $\lambda >0$ is a parameter. By using the minimax theorem
and the truncation technique, the existence of nontrivial solutions of Eq. (%
\ref{1-7}) is proved for $\lambda $ sufficiently large. Later, Liang and
Zhang \cite{LZ} obtained the existence and multiplicity of nontrivial
solutions for Eq. (\ref{1-7}) via Lions' second concentration-compactness
principle.

On the other hand, steep potential well has been applied to the study of the
existence and multiplicity of nontrivial solutions for biharmonic equations
without nonlocal term; see, for example, \cite{HLW,LCW,SCW,WZ,YT1}.
Specifically, Sun et. al. \cite{SCW} investigated the following biharmonic
equations with $p$-Laplacian and steep potential well
\begin{equation}
\left\{
\begin{array}{ll}
\Delta ^{2}u-\beta \Delta _{p}u+\lambda V(x)u=f(x,u) & \text{ in }\mathbb{R}%
^{N}, \\
u\in H^{2}(\mathbb{R}^{N}), &
\end{array}%
\right.  \label{1-8}
\end{equation}%
where $N\geq 1$, $\beta \in \mathbb{R}$, $\Delta _{p}u=\text{div}(|\nabla
u|^{p-2}\nabla u)$ with $p\geq 2$ and $\lambda V(x)$ is a steep potential
well. When $f$ satisfies various superlinear or sublinear assumptions, they
proved that Eq. (\ref{1-8}) admits one or two nontrivial solutions,
respectively.

Motivated by all results mentioned above, in the present paper we are
concerned with a class of generalized extensible beam equations with steep
potential well, i.e. Eq. $(E)$. We focus our attention on the multiplicity
and concentration of nontrivial solutions for Eq. $(E)$. Distinguished from
the existing literatures, (I) we allow the constant $b$ to be nonpositive,
which has the physical significance; (II) we are interested in seeking two
nontrivial solutions for Eq. $(E)$ with a superlinear nonlinearity, one of
which will blow up as the nonlocal term vanishes; (III) we would like to
explore the phenomenon of concentrations of two different nontrivial
solutions as $\lambda \rightarrow \infty $, which seems to be less involved
in extensible beam equations.

It is noteworthy that in analysis, we have to face some challenges. First,
since the constant $b\leq 0$ is allowed, how to construct an appropriate
norm of the working space such that this norm is associated with the norm $%
\Vert \nabla u\Vert _{L^{2}}=(\int_{\mathbb{R}^{N}}|\nabla u|^{2}dx)^{1/2}$
is crucial. Second, having considered the fact that the norms $\Vert \nabla
u\Vert _{L^{2}}$ and $\Vert u\Vert _{H^{2}}=(\int_{\mathbb{R}^{N}}(|\Delta
u|^{2}+|\nabla u|^{2}+u^{2})dx)^{1/2}$ are not equivalent, how to verify
that the energy functional of Eq. $(E)$ is bounded below and coercive in $%
H^{2}(\mathbb{R}^{N})$ is critical. Third, we note that $\Delta
u|_{\partial\Omega}=0$ is not included in the space $H_{0}^{1}(\Omega )\cap
H^{2}(\Omega ).$ In view of this, about the concentration of nontrivial
solutions, how to prove the functions of convergence satisfy the second
boundary condition $\Delta u|_{\partial\Omega}=0$ in Navier boundary
conditions is the key.

In order to overcome these difficulties, in this paper some new inequalities
are established and new research techniques are introduced. In addition, the
regularity of weak solutions for Navier BVPs to generalized extensible beam
equations is discussed. By so doing, we obtain the existence of two
nontrivial solutions for Eq. $(E)$ by the minimax theory and the
nonexistence of nontrivial solutions. Furthermore, we successfully figure
out the concentrations of two different nontrivial solutions for Eq. $(E)$
as $\lambda \rightarrow \infty $.

Before stating our results, we shall first introduce some notations. Denote
the best Sobolev constant for the imbedding $H^{2}(\mathbb{R}%
^{N})\hookrightarrow L^{r}(\mathbb{R}^{N})(2\leq r<+\infty )$ by $S_{r}$ for
$N=4.$ Let $A_{0}>0$ be a Gagliardo-Nirenberg constant satisfying the
following Gagliardo-Nirenberg inequality%
\begin{equation}
\int_{\mathbb{R}^{N}}\left\vert \nabla u\right\vert ^{2}dx\leq
A_{0}^{2}\left( \int_{\mathbb{R}^{N}}\left\vert \Delta u\right\vert
^{2}dx\right) ^{1/2}\left( \int_{\mathbb{R}^{N}}u^{2}dx\right) ^{1/2}.
\label{1-6}
\end{equation}%
Set%
\begin{equation*}
\beta _{N}:=\left\{
\begin{array}{ll}
\left( 1+\frac{A_{0}^{2}}{2}\right) \left( 1+\overline{A}_{N}^{16/N}\left%
\vert \left\{ a<c_{0}\right\} \right\vert ^{4/N}\right) & \text{ for }N=3,4,
\\
\left( 1+\frac{A_{0}^{2}}{2}\right) \left( 1+\overline{B}_{N}^{2}\left\vert
\left\{ V<c_{0}\right\} \right\vert ^{4/N}\right) & \text{ for }N\geq 4,%
\end{array}%
\right.
\end{equation*}%
and
\begin{equation*}
\Theta _{2,N}:=\left\{
\begin{array}{ll}
\left[ \left( 1+\frac{A_{0}^{2}}{2}\right) ^{-1}-S_{\infty }^{2}\left\vert
\left\{ V<c_{0}\right\} \right\vert \right] ^{-1} & \text{ for }N=3, \\
S_{2}^{-2}\left( 1+\frac{A_{0}^{2}}{2}\right) & \text{ for }N=4, \\
1+\frac{A_{0}^{2}}{2} & \text{ for }N>4.%
\end{array}%
\right.
\end{equation*}%
We now summarize our main results as follows.

\begin{theorem}
\label{t1-1}Suppose that $N\geq 3,\delta \geq \frac{2}{N-2}%
,b>-2A_{0}^{-2}\beta _{N}^{-1}$ and conditions $\left( V1\right) -\left(
V3\right) $ hold$.$ In addition, we assume that the function $f$ satisfies
the followings:

\begin{itemize}
\item[$(F1)$] $f(x,s)$ is a continuous function on $\mathbb{R}^{N}\times
\mathbb{R};$

\item[$(F2)$] there exists a constant $0<d_{0}<\alpha $ such that%
\begin{equation*}
pF(x,s)-f(x,s)s\leq d_{0}s^{2}\text{ for all }x\in \mathbb{R}^{N}\text{ and }%
s\in \mathbb{R},
\end{equation*}%
where%
\begin{equation*}
\alpha =\left\{
\begin{array}{ll}
\frac{1}{2}\delta \Theta _{2,N}^{-2}\left( 2+bA_{0}^{2}\beta _{N}\right) &
\text{ if }-2A_{0}^{-2}\beta _{N}^{-1}<b<0, \\
\delta \Theta _{2,N}^{-2} & \text{ if }b\geq 0,%
\end{array}%
\right.
\end{equation*}%
and $F(x,s)=\int_{0}^{s}f(x,t)dt;$

\item[$(F3)$] for each $\epsilon \in \left( 0,\frac{1}{2}\left(
2+bA_{0}^{2}\beta _{N}\right) \Theta _{2,N}^{-2}\right) ,$ there exist
constants $2<p<\frac{2N}{N-2}$ and $C_{1,\epsilon },C_{2,\epsilon }>0$
satisfying $C_{1,\epsilon }>\frac{2\delta +2-p}{\delta p}C_{2,\epsilon }$
such that for all $x\in \mathbb{R}^{N}\text{ and }s\in \mathbb{R}$,%
\begin{equation*}
C_{2,\epsilon }s^{p-1}-\gamma s\leq f\left( x,s\right) \leq \epsilon
s+C_{1,\epsilon }s^{p-1}
\end{equation*}%
for some constant $\gamma $ independent on $\epsilon .$
\end{itemize}

Then there exists constants $\Lambda _{1},a_{\ast }>0$ such that for every $%
\lambda \geq \Lambda _{1}$ and $0<a<a_{\ast },$ Eq. $(E)$ admits at least
two nontrivial solutions $u_{a,\lambda }^{(1)}$ and $u_{a,\lambda }^{(2)}$
satisfying $J_{a,\lambda }\left( u_{a,\lambda }^{(2)}\right) <0<J_{a,\lambda
}\left( u_{a,\lambda }^{(1)}\right) $. In particular, $u_{a,\lambda }^{(2)}$
is a ground state solution of Eq. $(E).$ Furthermore, when $\delta >\frac{2}{%
N-2},$ for every $\lambda \geq \Lambda _{1}$ there holds%
\begin{equation*}
J_{a,\lambda }\left( u_{a,\lambda }^{(2)}\right) \rightarrow -\infty \text{
and }\left\Vert u_{a,\lambda }^{(2)}\right\Vert _{\lambda }\rightarrow
\infty \text{ as }a\rightarrow 0,
\end{equation*}%
where $J_{a,\lambda }$ is the energy functional of Eq. $(E)$ and $\left\Vert
\cdot \right\Vert _{\lambda }$ is defined as (\ref{10}).
\end{theorem}

\begin{theorem}
\label{t1-2} Suppose that $N\geq 3,\delta \geq \frac{2}{N-2}%
,b>-2A_{0}^{-2}\beta _{N}^{-1}$ and conditions $(V1)-{(V2)}$ hold. In
addition, we assume that the function $f$ is a continuous function on $%
\mathbb{R}^{N}\times \mathbb{R}$ satisfying:

\begin{itemize}
\item[$(F3)^{\prime }$] for each $\epsilon \in \left( 0,b\overline{S}%
^{2}\left\vert \left\{ V<c_{0}\right\} \right\vert ^{-2/N}\right) ,$ there
exists constants $2<p<\frac{2N}{N-2}$ and $C_{1,\epsilon }>0$ such that for
all $x\in \mathbb{R}^{N}\text{ and }s\in \mathbb{R}$,%
\begin{equation*}
f\left( x,s\right) \leq \epsilon s+C_{1,\epsilon }s^{p-1}.
\end{equation*}
\end{itemize}

Then there exists $a^{\ast }>0$ such that for every $a>a^{\ast }$, Eq. $%
(K_{a,\lambda })$ does not admit any nontrivial solution for all $\lambda
>bc_{0}^{-1}\overline{S}^{2}\left\vert \left\{ V<c_{0}\right\} \right\vert
^{-2/N}.$
\end{theorem}

\begin{theorem}
\label{t1-3} Assume that $N\geq 5.$ Let $u_{a,\lambda }^{(1)}$ and $%
u_{a,\lambda }^{(2)}$ be the solutions obtained by Theorem \ref{t1-1}. Then $%
u_{a,\lambda }^{(1)}\rightarrow u_{\infty }^{(1)}$ and $u_{a,\lambda
}^{(2)}\rightarrow u_{\infty }^{(2)}$ in $H^{2}(\mathbb{R}^{N})$ as $\lambda
\rightarrow \infty ,$ where $u_{\infty }^{(1)}\neq u_{\infty }^{(2)}\in
H_{0}^{2}(\Omega )$ are two nontrivial solutions of the following Dirichlet
BVPs:%
\begin{equation}
\left\{
\begin{array}{ll}
\Delta ^{2}u-M\left( \int_{\Omega}\left\vert \nabla u\right\vert
^{2}dx\right) \Delta u=f(x,u) & \text{in }\Omega , \\
u=\frac{\partial u}{\partial n}=0, & \text{on}\ \partial \Omega .%
\end{array}%
\right.  \tag*{$\left( K_{\infty }\right) $}
\end{equation}
\end{theorem}

The remainder of this paper is organized as follows. After presenting some
preliminary results in section 2, we prove Theorem \ref{t1-1} in section 3,
and demonstrate proof of Theorem \ref{t1-2} in Sections 4. Sections 5 is
dedicated to the proof of Theorem \ref{t1-3}.

\section{Preliminaries}

Let%
\begin{equation*}
X=\left\{ u\in H^{2}(\mathbb{R}^{N})\ |\ \int_{\mathbb{R}^{N}}\left(
\left\vert \Delta u\right\vert ^{2}+V\left( x\right) u^{2}\right) dx<\infty
\right\}
\end{equation*}%
be equipped with the inner product and norm
\begin{equation*}
\left\langle u,v\right\rangle =\int_{\mathbb{R}^{N}}\left( \Delta u\Delta
v+V(x)uv\right) dx,\ \left\Vert u\right\Vert =\left\langle u,u\right\rangle
^{1/2}.
\end{equation*}%
For $\lambda >0$, we also need the following inner product and norm
\begin{equation}
\left\langle u,v\right\rangle _{\lambda }=\int_{\mathbb{R}^{N}}\left( \Delta
u\Delta v+\lambda V\left( x\right) uv\right) dx,\ \left\Vert u\right\Vert
_{\lambda }=\left\langle u,u\right\rangle _{\lambda }^{1/2}.  \label{10}
\end{equation}%
It is clear that $\left\Vert u\right\Vert \leq \left\Vert u\right\Vert
_{\lambda }$ for $\lambda \geq 1.$ Now we set $X_{\lambda }=(X,\left\Vert
u\right\Vert _{\lambda })$.

By the Young and Gagliardo-Nirenberg inequalities, there exists a sharp
constant $A_{0}>0$ such that%
\begin{equation}
\int_{\mathbb{R}^{N}}\left\vert \nabla u\right\vert ^{2}dx\leq \frac{%
A_{0}^{2}}{2}\int_{\mathbb{R}^{N}}\left( \left\vert \Delta u\right\vert
^{2}+u^{2}\right) dx.  \label{48}
\end{equation}%
This shows that
\begin{equation}
\int_{\mathbb{R}^{N}}\left( \left\vert \Delta u\right\vert ^{2}+u^{2}\right)
dx\leq \left\Vert u\right\Vert _{H^{2}}^{2}\leq \left( 1+\frac{A_{0}^{2}}{2}%
\right) \int_{\mathbb{R}^{N}}\left( \left\vert \Delta u\right\vert
^{2}+u^{2}\right) dx.  \label{1}
\end{equation}%
For $N=3,4,$ applying condition $(V1)$ and the H\"{o}lder, Young and
Gagliardo-Nirenberg inequalities, there exists a sharp constant $\overline{A}%
_{N}>0$ such that%
\begin{eqnarray}
\int_{\mathbb{R}^{N}}u^{2}dx &\leq &\frac{1}{c_{0}}\int_{\left\{ V\geq
c_{0}\right\} }V\left( x\right) u^{2}dx+\left( \left\vert \left\{
V<c_{0}\right\} \right\vert \int_{\mathbb{R}^{N}}\left\vert u\right\vert
^{4}dx\right) ^{\frac{1}{2}}  \notag \\
&\leq &\frac{1}{c_{0}}\int_{\mathbb{R}^{N}}V\left( x\right) u^{2}dx+\frac{N%
\overline{A}_{N}^{\frac{16}{N}}\left\vert \left\{ V<c_{0}\right\}
\right\vert ^{\frac{4}{N}}}{8}\int_{\mathbb{R}^{N}}\left\vert \Delta
u\right\vert ^{2}dx+\frac{8-N}{8}\int_{\mathbb{R}^{N}}u^{2}dx,  \notag
\end{eqnarray}%
which shows that%
\begin{equation}
\int_{\mathbb{R}^{N}}u^{2}dx\leq \frac{8}{Nc_{0}}\int_{\mathbb{R}%
^{N}}V\left( x\right) u^{2}dx+\overline{A}_{N}^{16/N}\left\vert \left\{
a<c_{0}\right\} \right\vert ^{4/N}\int_{\mathbb{R}^{N}}\left\vert \Delta
u\right\vert ^{2}dx.  \label{5}
\end{equation}%
It follows from (\ref{1}) and (\ref{5}) that
\begin{equation}
\left\Vert u\right\Vert _{H^{2}}^{2}\leq \left( 1+\frac{A_{0}^{2}}{2}\right)
\max \left\{ 1+\overline{A}_{N}^{16/N}\left\vert \left\{ a<c_{0}\right\}
\right\vert ^{4/N},\frac{8}{Nc_{0}}\right\} \left\Vert u\right\Vert ^{2}.
\label{40}
\end{equation}%
Similarly, we also obtain that%
\begin{equation}
\left\Vert u\right\Vert _{H^{2}}^{2}\leq \left( 1+\frac{A_{0}^{2}}{2}\right)
\left( 1+\overline{A}_{N}^{16/N}\left\vert \left\{ a<c_{0}\right\}
\right\vert ^{4/N}\right) \left\Vert u\right\Vert _{\lambda }^{2}  \label{44}
\end{equation}%
for $\lambda \geq 8N^{-1}c_{0}^{-1}\left( 1+\overline{A}_{N}^{16/N}\left%
\vert \left\{ a<c_{0}\right\} \right\vert ^{4/N}\right) .$ For $N>4,$ by
conditions $(V1)-(V2)$, H\"{o}lder and Gagliardo-Nirenberg inequalities,
there exists a sharp constant $\overline{B}_{N}>0$ such that
\begin{eqnarray*}
\int_{\mathbb{R}^{N}}u^{2}dx &=&\int_{\left\{ V\geq c_{0}\right\}
}u^{2}dx+\int_{\left\{ V<c_{0}\right\} }u^{2}dx \\
&\leq &\frac{1}{c_{0}}\int_{\mathbb{R}^{N}}V\left( x\right) u^{2}dx+%
\overline{B}_{N}^{2}\left\vert \left\{ V<c_{0}\right\} \right\vert
^{4/N}\int_{\mathbb{R}^{N}}\left\vert \Delta u\right\vert ^{2}dx.
\end{eqnarray*}%
Combining the above inequality with $\left( \ref{1}\right) $ yields
\begin{equation}
\left\Vert u\right\Vert _{H^{2}}^{2}\leq \left( 1+\frac{A_{0}^{2}}{2}\right)
\max \left\{ 1+\overline{B}_{N}^{2}\left\vert \left\{ V<c_{0}\right\}
\right\vert ^{4/N},\frac{1}{c_{0}}\right\} \left\Vert u\right\Vert ^{2}.
\label{41}
\end{equation}%
Similarly, we also have%
\begin{equation}
\left\Vert u\right\Vert _{H^{2}}^{2}\leq \left( 1+\frac{A_{0}^{2}}{2}\right)
\left( 1+\overline{B}_{N}^{2}\left\vert \left\{ V<c_{0}\right\} \right\vert
^{4/N}\right) \left\Vert u\right\Vert _{\lambda }^{2}  \label{45}
\end{equation}%
for $\lambda \geq c_{0}^{-1}\left( 1+\overline{B}_{N}^{2}\left\vert \left\{
V<c_{0}\right\} \right\vert ^{4/N}\right) .$ Set%
\begin{equation*}
\alpha _{N}:=\left\{
\begin{array}{ll}
\left( 1+\frac{A_{0}^{2}}{2}\right) \max \left\{ 1+\overline{A}%
_{N}^{16/N}\left\vert \left\{ a<c_{0}\right\} \right\vert ^{4/N},\frac{8}{%
Nc_{0}}\right\} & \text{ for }N=3,4, \\
\left( 1+\frac{A_{0}^{2}}{2}\right) \max \left\{ 1+\overline{B}%
_{N}^{2}\left\vert \left\{ V<c_{0}\right\} \right\vert ^{4/N},\frac{1}{c_{0}}%
\right\} & \text{ for }N\geq 5.%
\end{array}%
\right.
\end{equation*}%
Thus, it follows from (\ref{40}) and (\ref{41}) that
\begin{equation}
\left\Vert u\right\Vert _{H^{2}}^{2}\leq \alpha _{N}\left\Vert u\right\Vert
^{2},  \label{9}
\end{equation}%
which implies that the imbedding $X\hookrightarrow H^{2}(\mathbb{R}^{N})$ is
continuous. If we set%
\begin{equation*}
\Lambda _{N}:=\left\{
\begin{array}{ll}
8N^{-1}c_{0}^{-1}\left( 1+\overline{A}_{N}^{16/N}\left\vert \left\{
a<c_{0}\right\} \right\vert ^{4/N}\right) & \text{ for }N=3,4, \\
c_{0}^{-1}\left( 1+\overline{B}_{N}^{2}\left\vert \left\{ V<c_{0}\right\}
\right\vert ^{4/N}\right) & \text{ for }N\geq 5,%
\end{array}%
\right.
\end{equation*}%
then we have%
\begin{equation}
\left\Vert u\right\Vert _{H^{2}}^{2}\leq \beta _{N}\left\Vert u\right\Vert
_{\lambda }^{2}\text{ for }\lambda \geq \Lambda _{N},  \label{12}
\end{equation}%
where $\beta _{N}$ is defined as $(\ref{1-6})$. Furthermore, by $(\ref{48}),(%
\ref{1})$ and $(\ref{12})$ one has%
\begin{equation}
\int_{\mathbb{R}^{N}}\left\vert \nabla u\right\vert ^{2}dx\leq \frac{1}{2}%
A_{0}^{2}\beta _{N}\left\Vert u\right\Vert _{\lambda }^{2}\text{ for }%
\lambda \geq \Lambda _{N}.  \label{13}
\end{equation}%
Since the imbedding $H^{2}(\mathbb{R}^{3})\hookrightarrow L^{\infty }(%
\mathbb{R}^{3})$ is continuous, by $(\ref{44})$, for any $r\in \lbrack
2,+\infty )$ we have%
\begin{eqnarray}
\int_{\mathbb{R}^{3}}\left\vert u\right\vert ^{r}dx &\leq &\left\Vert
u\right\Vert _{L^{\infty }}^{r-2}\int_{\mathbb{R}^{3}}u^{2}dx  \notag \\
&\leq &S_{\infty }^{-\left( r-2\right) }\left( 1+\frac{A_{0}^{2}}{2}\right)
^{r/2}\left( 1+\overline{A}_{3}^{16/3}\left\vert \left\{ a<c_{0}\right\}
\right\vert ^{4/3}\right) ^{r/2}\left\Vert u\right\Vert _{\lambda }^{r}
\label{4}
\end{eqnarray}%
for $\lambda \geq \frac{8}{3c_{0}}\left( 1+\overline{A}_{3}^{16/3}\left\vert
\left\{ a<c_{0}\right\} \right\vert ^{4/3}\right) $. Moreover, using the
fact that the imbedding $H^{2}(\mathbb{R}^{4})\hookrightarrow L^{r}(\mathbb{R%
}^{4})$ $(2\leq r<+\infty )$ is continuous and $(\ref{44})$, for any $r\in
\lbrack 2,+\infty )$ one has
\begin{equation}
\int_{\mathbb{R}^{4}}\left\vert u\right\vert ^{r}dx\leq S_{r}^{-r}\left( 1+%
\frac{A_{0}^{2}}{2}\right) ^{r/2}\left( 1+\overline{A}_{4}^{4}\left\vert
\left\{ a<c_{0}\right\} \right\vert \right) ^{r/2}\left\Vert u\right\Vert
_{\lambda }^{r}  \label{6}
\end{equation}%
for $\lambda \geq 2c_{0}^{-1}\left( 1+\overline{A}_{4}^{4}\left\vert \left\{
a<c_{0}\right\} \right\vert \right) ,$ where $S_{r}$ is the best Sobolev
constant for the imbedding of $H^{2}(\mathbb{R}^{4})$ in $L^{r}(\mathbb{R}%
^{4})$ ($2\leq r<+\infty $). Finally, for $N>4,$ from conditions $(V1)-(V2)$%
, $(\ref{45})$ and H\"{o}lder and Gagliardo-Nirenberg inequalities again, it
follows that for any $r\in \lbrack 2,\frac{2N}{N-4}),$%
\begin{eqnarray}
&&\int_{\mathbb{R}^{N}}\left\vert u\right\vert ^{r}dx  \notag \\
&\leq &C_{0}^{N\left( r-2\right) /4}\left( \int_{\mathbb{R}^{N}}\left\vert
u\right\vert ^{2}dx\right) ^{[2N-r\left( N-4\right) ]/8}\left( \int_{\mathbb{%
R}^{N}}\left\vert \Delta u\right\vert ^{2}dx\right) ^{N\left( r-2\right) /8}
\notag \\
&\leq &C_{0}^{N\left( r-2\right) /4}\left( 1+\frac{A_{0}^{2}}{2}\right)
^{r/2}\left\Vert u\right\Vert _{\lambda }^{r}\text{ for }\lambda \geq \frac{%
1+C_{0}^{2}\left\vert \left\{ V<c_{0}\right\} \right\vert ^{4/N}}{c_{0}}.
\label{2}
\end{eqnarray}%
Set%
\begin{equation}
\Theta _{r,N}:=\left\{
\begin{array}{ll}
S_{\infty }^{-\left( r-2\right) }\left( 1+\frac{A_{0}^{2}}{2}\right)
^{r/2}\left( 1+\overline{A}_{3}^{16/3}\left\vert \left\{ a<c_{0}\right\}
\right\vert ^{4/3}\right) ^{r/2} & \text{ if }N=3, \\
S_{r}^{-r}\left( 1+\frac{A_{0}^{2}}{2}\right) ^{r/2}\left( 1+\overline{A}%
_{4}^{4}\left\vert \left\{ a<c_{0}\right\} \right\vert \right) ^{r/2} &
\text{ if }N=4, \\
C_{0}^{N\left( r-2\right) /4}\left( 1+\frac{A_{0}^{2}}{2}\right)
^{r/2}\left( 1+\overline{B}_{N}^{2}\left\vert \left\{ V<c_{0}\right\}
\right\vert ^{4/N}\right) ^{r/2} & \text{ if }N>4.%
\end{array}%
\right.  \label{7}
\end{equation}%
Thus, $\left( \ref{4}\right) -\left( \ref{7}\right) $ show that for any $%
r\in \lbrack 2,2_{\ast })$ and $\lambda \geq \Lambda _{N},$ there holds
\begin{equation}
\int_{\mathbb{R}^{N}}\left\vert u\right\vert ^{r}dx\leq \Theta
_{r,N}\left\Vert u\right\Vert _{\lambda }^{r}.  \label{11}
\end{equation}%
It is easily seen that Eq. $(E)$ is variational and its solutions are
critical points of the functional defined in $X_{\lambda }$ by
\begin{equation}
J_{a,\lambda }\left( u\right) =\frac{1}{2}\left\Vert u\right\Vert _{\lambda
}^{2}+\frac{a}{2\left( 1+\delta \right) }\Vert \nabla u\Vert
_{L^{2}}^{2(1+\delta )}+\frac{b}{2}\Vert \nabla u\Vert _{L^{2}}^{2}-\int_{%
\mathbb{R}^{N}}F\left( x,u\right) dx.  \label{2-3}
\end{equation}%
It is not difficult to prove that the functional $J_{a,\lambda }$ is of
class $C^{1}$ in $X_{\lambda }$, and that
\begin{eqnarray}
\langle J_{a,\lambda }^{\prime }(u),v\rangle &=&\int_{\mathbb{R}^{N}}\left[
\Delta u\cdot \Delta v+\lambda V\left( x\right) uv\right] dx+a\Vert \nabla
u\Vert _{L^{2}}^{2\delta }\int_{\mathbb{R}^{N}}\nabla u\cdot \nabla vdx
\notag \\
&&+b\int_{\mathbb{R}^{N}}\nabla u\cdot \nabla vdx-\int_{\mathbb{R}%
^{N}}f\left( x,u\right) vdx.  \label{2-2}
\end{eqnarray}%
Furthermore, we have the following results.

\begin{lemma}
\label{lem15}Suppose that $N\geq 3$ and $\delta \geq \frac{2}{N-2}.$ In
addition, we assume that conditions $(V1)-{(V2)},(F1)$ and ${(F3)}$ hold.
Then the energy functional $J_{a,\lambda }$ is bounded below and coercive on
$X_{\lambda }$ for all $a>0$ and%
\begin{equation*}
\lambda \geq \Lambda _{0}:=\left\{
\begin{array}{ll}
\max \left\{ \Lambda _{N},\frac{2\epsilon }{c_{0}}\right\} & \text{ if }%
\delta >\frac{2}{N-2}, \\
\max \left\{ \Lambda _{N},\frac{2\epsilon }{c_{0}}+\frac{4C_{1,\epsilon }}{%
c_{0}p}\left( \frac{2C_{1,\epsilon }\left( 1+\delta \right) }{ap\overline{S}%
^{2N/(N-2)}}\right) ^{\frac{(p-2)(N-2)}{2N-p\left( N-2\right) }}\right\} &
\text{ if }\delta =\frac{2}{N-2}.%
\end{array}%
\right.
\end{equation*}%
Furthermore, for all $a>0$ and $\lambda \geq \Lambda _{0},$ there exists a
constant $R_{a}>0$ such that
\begin{equation*}
J_{a,\lambda }(u)\geq 0\text{ for all }u\in X_{\lambda }\text{ with }%
\left\Vert u\right\Vert _{\lambda }\geq R_{a}
\end{equation*}
\end{lemma}

\textbf{Proof.} Let $u\in X_{\lambda }.$ Note that for any $2\leq r\leq
2^{\ast }:=\frac{2N}{N-2},$ there holds%
\begin{eqnarray}
&&\int_{\mathbb{R}^{N}}\left\vert u\right\vert ^{r}dx  \notag \\
&\leq &\left( \int_{\mathbb{R}^{N}}\left\vert u\right\vert ^{2}dx\right) ^{%
\frac{2^{\ast }-r}{2^{\ast }-2}}\left( \int_{\mathbb{R}^{N}}\left\vert
u\right\vert ^{2^{\ast }}dx\right) ^{\frac{r-2}{2^{\ast }-2}}  \notag \\
&\leq &\left( \frac{1}{\lambda c_{0}}\int_{\mathbb{R}^{N}}\lambda
V(x)u^{2}dx+\overline{S}^{-2}\left\vert \left\{ V<c_{0}\right\} \right\vert
^{\frac{2}{N}}\Vert \nabla u\Vert _{L^{2}}^{2}\right) ^{\frac{2^{\ast }-r}{%
2^{\ast }-2}}\left( \overline{S}^{-1}\Vert \nabla u\Vert _{L^{2}}\right) ^{%
\frac{N\left( r-2\right) }{2}},  \label{8}
\end{eqnarray}%
where we have used the H\"{o}lder and Sobolev inequalities and $\overline{S}$
is the best Sobolev constant for the imbedding of $D^{1,2}(R^{N})$ in $%
L^{2^{\ast }}(R^{N})$. We now divide the proof into two separate cases:%
\newline
Case $A:\int_{\mathbb{R}^{N}}\lambda V(x)u^{2}dx\geq \lambda c_{0}\left(
\frac{4C_{1,\epsilon }}{p(\lambda c_{0}-2\epsilon )}\right) ^{\frac{4}{%
(p-2)(N-2)}}\left( \overline{S}^{-1}\Vert \nabla u\Vert _{L^{2}}\right)
^{2^{\ast }}.$ It follows from condition $(F3)$ and $(\ref{8})$ that%
\begin{eqnarray*}
&&J_{a,\lambda }\left( u\right) \\
&\geq &\frac{1}{2}\left\Vert u\right\Vert _{\lambda }^{2}+\frac{a}{2\left(
1+\delta \right) }\Vert \nabla u\Vert _{L^{2}}^{2(1+\delta )}+\frac{b}{2}%
\Vert \nabla u\Vert _{L^{2}}^{2}-\frac{\epsilon }{2}\int_{\mathbb{R}%
^{N}}u^{2}dx-\frac{C_{1,\epsilon }}{p}\int_{\mathbb{R}^{N}}\left\vert
u\right\vert ^{p}dx \\
&\geq &\frac{1}{4}\left\Vert u\right\Vert _{\lambda }^{2}+\frac{a}{2\left(
1+\delta \right) }\Vert \nabla u\Vert _{L^{2}}^{2(1+\delta )}+\frac{1}{2}%
\left( b-\epsilon \overline{S}^{-2}\left\vert \left\{ V<c_{0}\right\}
\right\vert ^{\frac{2}{N}}\right) \Vert \nabla u\Vert _{L^{2}}^{2} \\
&&-\frac{C_{1,\epsilon }}{p\overline{S}^{p}}\left\vert \left\{
V<c_{0}\right\} \right\vert ^{\frac{2N-p\left( N-2\right) }{2N}}\Vert \nabla
u\Vert _{L^{2}}^{p}.
\end{eqnarray*}%
Since $\delta \geq \frac{2}{N-2},$ we have $1+\delta >\frac{p}{2}>1.$ Then
there exists a constant $D_{a}$ such that%
\begin{eqnarray*}
D_{a} &=&\min_{t\geq 0}\left[ \frac{at^{1+\delta }}{2\left( 1+\delta \right)
}+\frac{t}{2}\left( b-\frac{\epsilon \left\vert \left\{ V<c_{0}\right\}
\right\vert ^{\frac{2}{N}}}{\overline{S}^{2}}\right) -\frac{C_{1,\epsilon
}t^{\frac{p}{2}}}{p\overline{S}^{p}}\left\vert \left\{ V<c_{0}\right\}
\right\vert ^{\frac{2N-p\left( N-2\right) }{2N}}\right] \\
&<&0,
\end{eqnarray*}%
and $D_{a}\rightarrow -\infty $ as $a\rightarrow 0$. Using this, together
with the above inequality leads to%
\begin{equation*}
J_{a,\lambda }\left( u\right) \geq \frac{1}{4}\left\Vert u\right\Vert
_{\lambda }^{2}+D_{a}\geq D_{a},
\end{equation*}%
which implies that $J_{a,\lambda }\left( u\right) $ is bounded below and
coercive on $X_{\lambda }$ for all $a>0$ and $\lambda >\max \left\{ \Lambda
_{N},\frac{2\epsilon }{c_{0}}\right\} .$\newline
Case $B:\int_{\mathbb{R}^{N}}\lambda V(x)u^{2}dx<\lambda c_{0}\left( \frac{%
4C_{1,\epsilon }}{p(\lambda c_{0}-2\epsilon )}\right) ^{\frac{4}{(p-2)(N-2)}%
}\left( \overline{S}^{-1}\Vert \nabla u\Vert _{L^{2}}\right) ^{2^{\ast }}.$
By virtue of $(\ref{8})$ one has%
\begin{eqnarray*}
&&\int_{\mathbb{R}^{N}}|u|^{p}dx \\
&\leq &\left( \frac{1}{\lambda c_{0}}\int_{\mathbb{R}^{N}}\lambda
V(x)u^{2}dx+\frac{\left\vert \left\{ V<c_{0}\right\} \right\vert ^{\frac{2}{N%
}}}{\overline{S}^{2}}\Vert \nabla u\Vert _{L^{2}}^{2}\right) ^{\frac{2^{\ast
}-p}{2^{\ast }-2}}\cdot \left( \overline{S}^{-1}\Vert \nabla u\Vert
_{L^{2}}\right) ^{\frac{N\left( p-2\right) }{2}} \\
&\leq &\overline{S}^{-2^{\ast }}\left( \frac{4C_{1,\epsilon }}{p(\lambda
c_{0}-2\epsilon )}\right) ^{\frac{2N-p\left( N-2\right) }{(p-2)(N-2)}}\Vert
\nabla u\Vert _{L^{2}}^{2^{\ast }}+\frac{\left\vert \left\{ V<c_{0}\right\}
\right\vert ^{\frac{2N-p\left( N-2\right) }{2N}}}{\overline{S}^{p}}\Vert
\nabla u\Vert _{L^{2}}^{p}.
\end{eqnarray*}%
Using this, together with condition $(F3),$ gives%
\begin{eqnarray*}
J_{a,\lambda }\left( u\right) &\geq &\frac{1}{2}\left\Vert u\right\Vert
_{\lambda }^{2}+\frac{a}{2\left( 1+\delta \right) }\Vert \nabla u\Vert
_{L^{2}}^{2(1+\delta )}+\frac{b}{2}\Vert \nabla u\Vert _{L^{2}}^{2} \\
&&-\frac{\epsilon }{2}\int_{\mathbb{R}^{N}}u^{2}dx-\frac{C_{1,\epsilon }}{p}%
\int_{\mathbb{R}^{N}}\left\vert u\right\vert ^{p}dx \\
&\geq &\frac{1}{4}\left\Vert u\right\Vert _{\lambda }^{2}+\frac{a}{2\left(
1+\delta \right) }\Vert \nabla u\Vert _{L^{2}}^{2(1+\delta )}+\frac{1}{2}%
\left( b-\epsilon \overline{S}^{-2}\left\vert \left\{ V<c_{0}\right\}
\right\vert ^{\frac{2}{N}}\right) \Vert \nabla u\Vert _{L^{2}}^{2} \\
&&-\frac{C_{1,\epsilon }}{p\overline{S}^{2^{\ast }}}\left( \frac{%
4C_{1,\epsilon }}{p(\lambda c_{0}-2\epsilon )}\right) ^{\frac{2N-p\left(
N-2\right) }{(p-2)(N-2)}}\Vert \nabla u\Vert _{L^{2}}^{2^{\ast }}-\frac{%
C_{1,\epsilon }}{p\overline{S}^{p}}\left\vert \left\{ V<c_{0}\right\}
\right\vert ^{\frac{2N-p\left( N-2\right) }{2N}}\Vert \nabla u\Vert
_{L^{2}}^{p}.
\end{eqnarray*}%
If $\delta =\frac{2}{N-2},$ then for%
\begin{equation*}
\lambda >\frac{2\epsilon }{c_{0}}+\frac{4C_{1,\epsilon }}{c_{0}p}\left[
\frac{2C_{1,\epsilon }\left( 1+\delta \right) }{ap\overline{S}^{2^{\ast }}}%
\right] ^{\frac{(p-2)(N-2)}{2N-p\left( N-2\right) }},
\end{equation*}%
there exists a constant $\overline{D}_{a}<D_{a}<0$ such that%
\begin{eqnarray*}
J_{a,\lambda }\left( u\right) &\geq &\frac{1}{4}\left\Vert u\right\Vert
_{\lambda }^{2}+\frac{1}{2}\left( b-\epsilon \overline{S}^{-2}\left\vert
\left\{ V<c_{0}\right\} \right\vert ^{\frac{2}{N}}\right) \Vert \nabla
u\Vert _{L^{2}}^{2} \\
&&+\left[ \frac{a}{2\left( 1+\delta \right) }-\frac{C_{1,\epsilon }}{p%
\overline{S}^{2^{\ast }}}\left( \frac{4C_{1,\epsilon }}{p(\lambda
c_{0}-2\epsilon )}\right) ^{\frac{2N-p\left( N-2\right) }{(p-2)(N-2)}}\right]
\Vert \nabla u\Vert _{L^{2}}^{2(1+\delta )} \\
&&-\frac{C_{1,\epsilon }}{p\overline{S}^{p}}\left\vert \left\{
V<c_{0}\right\} \right\vert ^{\frac{2N-p\left( N-2\right) }{2N}}\Vert \nabla
u\Vert _{L^{2}}^{p} \\
&\geq &\frac{1}{4}\left\Vert u\right\Vert _{\lambda }^{2}+\overline{D}%
_{a}\geq \overline{D}_{a}.
\end{eqnarray*}%
If $\delta >\frac{2}{N-2},$ then for $\lambda >\frac{2\epsilon }{c_{0}},$
there exists a constant $\widetilde{D}_{a}<0$ such that%
\begin{eqnarray*}
J_{a,\lambda }\left( u\right) &\geq &\frac{1}{4}\left\Vert u\right\Vert
_{\lambda }^{2}+\frac{a}{2\left( 1+\delta \right) }\Vert \nabla u\Vert
_{L^{2}}^{2(1+\delta )}+\frac{1}{2}\left( b-\epsilon \overline{S}%
^{-2}\left\vert \left\{ V<c_{0}\right\} \right\vert ^{\frac{2}{N}}\right)
\Vert \nabla u\Vert _{L^{2}}^{2} \\
&&-\frac{C_{1,\epsilon }}{p\overline{S}^{2^{\ast }}}\left( \frac{%
4C_{1,\epsilon }}{p(\lambda c_{0}-2\epsilon )}\right) ^{\frac{2N-p\left(
N-2\right) }{(p-2)(N-2)}}\Vert \nabla u\Vert _{L^{2}}^{2^{\ast }}-\frac{%
C_{1,\epsilon }}{p\overline{S}^{p}}\left\vert \left\{ V<c_{0}\right\}
\right\vert ^{\frac{2N-p\left( N-2\right) }{2N}}\Vert \nabla u\Vert
_{L^{2}}^{p} \\
&\geq &\widetilde{D}_{a}.
\end{eqnarray*}%
This indicates that $J_{a,\lambda }$ is bounded below and coercive on $%
X_{\lambda }$ for all $a>0$ and $\lambda \geq \Lambda _{0}.$ Furthermore,
for all $a>0$ and $\lambda \geq \Lambda _{0},$ it is clear that there exists
a constant $R_{a}>0$ such that
\begin{equation*}
J_{a,\lambda }(u)\geq 0\text{ for all }u\in X_{\lambda }\text{ with }%
\left\Vert u\right\Vert _{\lambda }\geq R_{a}.
\end{equation*}%
Consequently, the proof is complete. \hfill $\square $

Next, we give a useful theorem, which is the variant version of the mountain
pass theorem. It can help us to find a so-called Cerami type $(PS)$ sequence.

\begin{lemma}
\label{l2}(\cite{E1}, Mountain Pass Theorem). Let $E$ be a real Banach space
with its dual space $E^{\ast },$ and suppose that $I\in C^{1}(E,R)$
satisfies
\begin{equation*}
\max \{I(0),I(e)\}\leq \mu <\eta \leq \inf_{\Vert u\Vert =\rho }I(u),
\end{equation*}%
for some $\mu <\eta ,\rho >0$ and $e\in E$ with $\Vert e\Vert >\rho $. Let $%
c\geq \eta $ be characterized by
\begin{equation*}
c=\inf_{\gamma \in \Gamma }\max_{0\leq \tau \leq 1}I(\gamma (\tau )),
\end{equation*}%
where $\Gamma =\{\gamma \in C([0,1],E):\gamma (0)=0,\gamma (1)=e\}$ is the
set of continuous paths joining $0$ and $e$, then there exists a sequence $%
\{u_{n}\}\subset E$ such that
\begin{equation*}
I(u_{n})\rightarrow c\geq \eta \quad \text{and}\quad (1+\Vert u_{n}\Vert
)\Vert I^{\prime }(u_{n})\Vert _{E^{\ast }}\rightarrow 0\quad \text{as}\
n\rightarrow \infty .
\end{equation*}
\end{lemma}

In what follows, we give two lemmas which ensure that the functional $%
J_{a,\lambda }$ has the mountain pass geometry.

\begin{lemma}
\label{lem1}Suppose that $b>-2A_{0}^{-2}\beta _{N}^{-1}.$ In addition,
assume that conditions $(V1)-(V2),(F1)$ and $(F3)$ hold. Then there exists $%
\rho >0$ such that for every $a>0$ and $\lambda >\Lambda _{N},$%
\begin{equation*}
\inf \{J_{a,\lambda }(u):u\in X_{\lambda }\ \text{with}\ \Vert u\Vert
_{\lambda }=\rho \}>\eta
\end{equation*}%
for some $\eta >0.$
\end{lemma}

\noindent Proof. \/ By $(\ref{13})$ and the condition $(F3),$ for all $u\in
X_{\lambda }$ one has%
\begin{eqnarray*}
J_{a,\lambda }(u) &\geq &\frac{1}{2}\left\Vert u\right\Vert _{\lambda }^{2}+%
\frac{a}{2\left( 1+\delta \right) }\Vert \nabla u\Vert _{L^{2}}^{2(1+\delta
)}+\frac{b}{2}\Vert \nabla u\Vert _{L^{2}}^{2}-\frac{\epsilon }{2}\int_{%
\mathbb{R}^{N}}u^{2}dx-\frac{C_{1,\epsilon }}{p}\int_{\mathbb{R}%
^{N}}\left\vert u\right\vert ^{p}dx \\
&\geq &\left\{
\begin{array}{ll}
\frac{1}{2}\left( 1-\epsilon \Theta _{2,N}^{2}\right) \Vert u\Vert _{\lambda
}^{2}-\frac{C_{1,\epsilon }}{p}\Theta _{p,N}^{p}\left\Vert u\right\Vert
_{\lambda }^{p} & \text{ if }b\geq 0, \\
\frac{1}{2}\left( 1+\frac{bA_{0}^{2}}{2}\beta _{N}-\epsilon \Theta
_{2,N}^{2}\right) \Vert u\Vert _{\lambda }^{2}-\frac{C_{1,\epsilon }}{p}%
\Theta _{p,N}^{p}\left\Vert u\right\Vert _{\lambda }^{p} & \text{ if }%
-2A_{0}^{-2}\beta _{N}^{-1}<b<0.%
\end{array}%
\right.
\end{eqnarray*}%
Let%
\begin{equation*}
g(t)=\frac{1}{2}\left( 1-\epsilon \Theta _{2,N}^{2}\right) t^{2}-\frac{%
C_{1,\epsilon }\Theta _{p,N}^{p}}{p}t^{p}\text{ for }t\geq 0.
\end{equation*}%
A direct calculation shows that%
\begin{equation*}
\max_{t\geq 0}g(t)=g(\bar{t})=\frac{(p-2)}{2p}\left( 1-\epsilon \Theta
_{2,N}^{2}\right) ^{p/(p-2)}\left( C_{1,\epsilon }\Theta _{p,N}^{p}\right)
^{-2/(p-2)},
\end{equation*}%
where
\begin{equation*}
\bar{t}=\left[ \frac{\left( 1-\epsilon \Theta _{2,N}^{2}\right) }{%
C_{1,\epsilon }\Theta _{p,N}^{p}}\right] ^{1/(p-2)}.
\end{equation*}%
This shows that when $b\geq 0,$ for every $u\in X_{\lambda }$ with $%
\left\Vert u\right\Vert _{\lambda }=\bar{t}\ $we have%
\begin{equation*}
J_{a,\lambda }\left( u\right) \geq g(\bar{t})>0.
\end{equation*}%
Choosing $\rho =\bar{t}$ and%
\begin{equation*}
\eta =\frac{(p-2)}{2p}\left( 1-\epsilon \Theta _{2,N}^{2}\right)
^{p/(p-2)}\left( C_{1,\epsilon }\Theta _{p,N}^{p}\right) ^{-2/(p-2)}>0,
\end{equation*}%
it is easy to see that the result holds. Similarly, when $-2A_{0}^{-2}\beta
_{N}^{-1}<b<0,$ for every $u\in X_{\lambda }$ with%
\begin{equation*}
\left\Vert u\right\Vert _{\lambda }=\tilde{t}=\left[ \frac{\left( 1+\frac{%
bA_{0}^{2}}{2}\beta _{N}-\epsilon \Theta _{2,N}^{2}\right) }{C_{1,\epsilon
}\Theta _{p,N}^{p}}\right] ^{1/(p-2)},
\end{equation*}%
we can take $\rho =\tilde{t}$ and
\begin{equation*}
\eta =\frac{(p-2)}{2p}\left( 1+\frac{bA_{0}^{2}}{2}\beta _{N}-\epsilon
\Theta _{2,N}^{2}\right) ^{p/(p-2)}\left( C_{1,\epsilon }\Theta
_{p,N}^{p}\right) ^{-2/(p-2)}
\end{equation*}%
such that the result holds. This completes the proof. \hfill $\square $

Define%
\begin{equation}
\Pi _{\lambda }=\sup_{u\in X_{\lambda }\backslash \{0\}}\frac{\left( \int_{%
\mathbb{R}^{N}}\left\vert u\right\vert ^{p}dx\right) ^{1/p}}{\left\Vert
u\right\Vert _{\lambda }}.  \label{28}
\end{equation}%
It follows from $\left( \ref{11}\right) $ that%
\begin{equation}
\Pi _{\lambda }\leq \Theta _{p,N}\text{ for }\lambda \geq \Lambda _{N}.
\label{29}
\end{equation}%
Furthermore, by Appendix A there exist $\Lambda _{1}\geq \Lambda _{N}$ and $%
\phi _{\lambda }\in X_{\lambda }\backslash \{0\}$ such that
\begin{equation}
\Pi _{\lambda }=\frac{\left( \int_{\mathbb{R}^{N}}\left\vert \phi _{\lambda
}\right\vert ^{p}dx\right) ^{1/p}}{\left\Vert \phi _{\lambda }\right\Vert
_{\lambda }}>0\text{ for every }\lambda \geq \Lambda _{1},  \label{31}
\end{equation}%
and there exists a constant $\Pi _{\infty }>0$ independent on $\lambda $
such that%
\begin{equation}
\Pi _{\lambda }\searrow \Pi _{\infty }\text{ as }\lambda \nearrow \infty .
\label{32}
\end{equation}%
Setting
\begin{equation*}
a_{\ast }:=\frac{2^{2+\delta }C_{2,\epsilon }\Pi _{\infty }^{p}\left(
1+\delta \right) (p-2)}{\delta pA_{0}^{2(1+\delta )}\beta _{N}^{1+\delta }}%
\left[ \frac{C_{2,\epsilon }\Pi _{\infty }^{p}(2\delta +2-p)}{\delta p\left(
1+\frac{bA_{0}^{2}}{2}\beta _{N}+\gamma \Theta _{2,N}^{2}\right) }\right] ^{%
\frac{2\delta +2-p}{p-2}}.
\end{equation*}

\begin{lemma}
\label{lem2}Assume that $b\in R,$ conditions $(V1)-(V3),(F1)$ and $(F3)$
hold. Let $\rho >0$ be as in Lemma \ref{lem1}. Then for every $\lambda \geq
\Lambda _{1}$ and $0<a<a_{\ast },$ there exists $e\in X_{\lambda }$
satisfying
\begin{equation*}
\Vert e\Vert _{\lambda }>\rho \ \text{and}\ \Vert e\Vert _{\lambda
}\rightarrow \infty \ \text{as}\ a\rightarrow 0
\end{equation*}%
such that%
\begin{equation*}
J_{a,\lambda }(e)<0\ \text{and}\ J_{a,\lambda }(e)\rightarrow -\infty \
\text{as}\ a\rightarrow 0.
\end{equation*}
\end{lemma}

\noindent Proof. \/ Let $\phi _{\lambda }\in X_{\lambda }\backslash \{0\}$
be as in $\left( \ref{31}\right) $ and let
\begin{eqnarray*}
I\left( t\right) &=&I_{a,\lambda }\left( t\phi _{\lambda }\right) \\
&=&\frac{t^{2}}{2}\left\Vert \phi _{\lambda }\right\Vert _{\lambda }^{2}+%
\frac{at^{2(1+\delta )}}{2\left( 1+\delta \right) }\Vert \nabla \phi
_{\lambda }\Vert _{L^{2}}^{2(1+\delta )}+\frac{bt^{2}}{2}\Vert \nabla \phi
_{\lambda }\Vert _{L^{2}}^{2} \\
&&+\frac{\gamma t^{2}}{2}\int_{\mathbb{R}^{N}}\phi _{\lambda }^{2}dx-\frac{%
C_{2,\epsilon }t^{p}}{p}\int_{\mathbb{R}^{N}}\left\vert \phi _{\lambda
}\right\vert ^{p}dx\text{ for }t>0.
\end{eqnarray*}%
Then it follows from (\ref{13}) and (\ref{11}) that%
\begin{eqnarray*}
I\left( t\right) &\leq &\frac{A_{0}^{2(1+\delta )}\beta _{N}^{1+\delta
}\left\Vert \phi _{\lambda }\right\Vert _{\lambda }^{2(1+\delta )}t^{2}}{%
2^{3+\delta }\left( 1+\delta \right) } \\
&&\cdot \left[ at^{2\delta }+\frac{2^{2+\delta }\left( 1+\delta \right)
\left( 1+\frac{bA_{0}^{2}}{2}\beta _{N}+\gamma \Theta _{2,N}^{2}\right) }{%
A_{0}^{2(1+\delta )}\beta _{N}^{1+\delta }\left\Vert \phi _{\lambda
}\right\Vert _{\lambda }^{2\delta }}-\frac{2^{3+\delta }\left( 1+\delta
\right) C_{2,\epsilon }\Pi _{\infty }^{p}}{pA_{0}^{2(1+\delta )}\beta
_{N}^{1+\delta }\left\Vert \phi _{\lambda }\right\Vert _{\lambda
}^{2(1+\delta )-p}}t^{p-2}\right] .
\end{eqnarray*}%
A direct calculation shows that there exists%
\begin{equation*}
t_{a,\lambda }:=\left( \frac{2^{2+\delta }C_{2,\epsilon }\Pi _{\infty
}^{p}\left( 1+\delta \right) (p-2)}{a\delta pA_{0}^{2(1+\delta )}\beta
_{N}^{1+\delta }}\right) ^{1/\left( 2\delta +2-p\right) }\left\Vert \phi
_{\lambda }\right\Vert _{\lambda }^{-1}>0
\end{equation*}%
such that for every $0<a<a_{\ast },$
\begin{eqnarray*}
&&at_{a,\lambda }^{2\delta }+\frac{2\left( 1+\delta \right) \left( b+\gamma
\Theta _{2,N}^{2}\right) }{A_{0}^{2}\beta _{N}\left\Vert \phi _{\lambda
}\right\Vert _{\lambda }^{2\delta }}-\frac{4\left( 1+\delta \right)
C_{2,\epsilon }\Pi _{\infty }^{p}}{pA_{0}^{2}\beta _{N}\left\Vert \phi
_{\lambda }\right\Vert _{\lambda }^{2(1+\delta )-p}}t_{a,\lambda }^{p-2} \\
&=&\frac{2^{2+\delta }\left( 1+\delta \right) \left\Vert \phi _{\lambda
}\right\Vert _{\lambda }^{-2\delta }}{A_{0}^{2(1+\delta )}\beta
_{N}^{1+\delta }} \\
&&\cdot \left[ \left( 1+\frac{bA_{0}^{2}}{2}\beta _{N}+\gamma \Theta
_{2,N}^{2}\right) -\frac{C_{2,\epsilon }\Pi _{\infty }^{p}(2\delta +2-p)}{%
\delta p}\left( \frac{2^{2+\delta }C_{2,\epsilon }\Pi _{\infty }^{p}\left(
1+\delta \right) (p-2)}{a\delta pA_{0}^{2(1+\delta )}\beta _{N}^{1+\delta }}%
\right) ^{\frac{p-2}{2\delta +2-p}}\right] \\
&<&0,
\end{eqnarray*}%
this implies that%
\begin{equation*}
I\left( t_{a,\lambda }\right) =I_{a,\lambda }\left( t_{a,\lambda }\phi
_{\lambda }\right) <0\text{ for }0<a<a_{\ast }
\end{equation*}%
and
\begin{equation*}
I_{a,\lambda }\left( t_{a,\lambda }\phi _{\lambda }\right) \rightarrow
-\infty \text{ as }a\rightarrow 0.
\end{equation*}%
Choosing $e=t_{a,\lambda }|\phi _{\lambda }|$. Clearly,
\begin{equation*}
\Vert e\Vert _{\lambda }=\Vert t_{a,\lambda }\phi _{\lambda }\Vert _{\lambda
}=\left[ \frac{2^{2+\delta }C_{2,\epsilon }\Pi _{\infty }^{p}\left( 1+\delta
\right) (p-2)}{a\delta pA_{0}^{2(1+\delta )}\beta _{N}^{1+\delta }}\right]
^{1/\left( 2\delta +2-p\right) }\rightarrow \infty \ \text{as}\ a\rightarrow
0.
\end{equation*}%
Note that for $0<a<a_{\ast },$ by $\left( \ref{29}\right) $ and $\left( \ref%
{32}\right) ,$ there holds%
\begin{equation*}
\left[ \frac{2^{2+\delta }C_{2,\epsilon }\Pi _{\infty }^{p}\left( 1+\delta
\right) (p-2)}{a\delta pA_{0}^{2(1+\delta )}\beta _{N}^{1+\delta }}\right] ^{%
\frac{1}{2\delta +2-p}}>\left[ \frac{\delta p\left( 1+\frac{bA_{0}^{2}}{2}%
\beta _{N}+\gamma \Theta _{2,N}^{2}\right) }{C_{2,\epsilon }\Theta
_{p,N}^{p}(2\delta +2-p)}\right] ^{\frac{1}{p-2}},
\end{equation*}%
by using $(\ref{29}).$ Using this, together with condition $(F3),$ leads to%
\begin{equation*}
\Vert e\Vert _{\lambda }>\rho :=\left\{
\begin{array}{ll}
\left( \frac{1-\epsilon \Theta _{2,N}^{2}}{C_{1,\epsilon }\Theta _{p,N}^{p}}%
\right) ^{1/(p-2)} & \text{ if }b\geq 0, \\
\left( \frac{1+\frac{bA_{0}^{2}}{2}\beta _{N}-\epsilon \Theta _{2,N}^{2}}{%
C_{1,\epsilon }\Theta _{p,N}^{p}}\right) ^{\frac{1}{p-2}} & \text{ if }%
-2A_{0}^{-2}\beta _{N}^{-1}<b<0%
\end{array}%
\right. ,
\end{equation*}%
where $\rho >0$ is as in Lemma \ref{lem1}. Moreover, by condition $(F3)$,
there holds $J_{a,\lambda }(e)\leq I_{a,\lambda }(e)<0$ for $0<a<a_{\ast }$.
Consequently, the lemma is proved. \hfill $\square $

Define%
\begin{equation*}
c_{\lambda }=\inf_{\gamma \in \Gamma _{\lambda }}\max_{0\leq t\leq
1}J_{a,\lambda }(\gamma (t))
\end{equation*}%
and%
\begin{equation*}
c_{0}(\Omega )=\inf_{\gamma \in \overline{\Gamma }_{\lambda }(\Omega
)}\max_{0\leq t\leq 1}J_{a,\lambda }|_{H_{0}^{2}(\Omega )}(\gamma (t)),
\end{equation*}%
where $J_{a,\lambda }|_{H_{0}^{2}(\Omega )}$ is a restriction of $%
J_{a,\lambda }$ on $H_{0}^{2}(\Omega ),$%
\begin{equation*}
\Gamma _{\lambda }=\{\gamma \in C([0,1],X_{\lambda }):\gamma (0)=0,\gamma
(1)=e\}
\end{equation*}%
and%
\begin{equation*}
\overline{\Gamma }_{\lambda }(\Omega )=\{\gamma \in C([0,1],H_{0}^{2}(\Omega
)):\gamma (0)=0,\gamma (1)=e\}.
\end{equation*}%
Note that for $u\in H_{0}^{2}(\Omega ),$
\begin{equation*}
J_{a,\lambda }|_{H_{0}^{2}(\Omega )}(u)=\frac{1}{2}\int_{\Omega }\left\vert
\Delta u\right\vert ^{2}dx+\frac{a}{2\left( 1+\delta \right) }\left(
\int_{\Omega }\left\vert \nabla u\right\vert ^{2}dx\right) ^{2(1+\delta )}+%
\frac{b}{2}\int_{\Omega }\left\vert \nabla u\right\vert ^{2}dx-\int_{\Omega
}F(x,u)dx
\end{equation*}%
and $c_{0}(\Omega )$ independent of $\lambda .$ Moreover, if conditions $(F1)
$ and $(F3)$ hold, then by the proofs of Lemmas \ref{lem1} and \ref{lem2},
we can conclude that $J_{a,\lambda }|_{H_{0}^{2}(\Omega )}$ satisfies the
mountain pass hypothesis as in Theorem \ref{l2}.

Since $H_{0}^{2}(\Omega )\subset X_{\lambda }$ for all $\lambda >0$, one can
see that $0<\eta \leq c_{\lambda }\leq c_{0}(\Omega )$ for all $\lambda \geq
\Lambda _{N}.$ Take $D_{0}>c_{0}(\Omega ).$ Then we have
\begin{equation*}
0<\eta \leq c_{\lambda }\leq c_{0}(\Omega )<D_{0}\text{ for all }\lambda
\geq \Lambda _{N}.
\end{equation*}%
By Lemmas \ref{lem1}, \ref{lem2} and Theorem \ref{l2}, we obtain that for
each $\lambda \geq \Lambda _{N}$, there exists a sequence $\{u_{n}\}\subset
X_{\lambda }$ such that
\begin{equation*}
J_{a,\lambda }(u_{n})\rightarrow c_{\lambda }>0\quad \text{and}\quad
(1+\Vert u_{n}\Vert _{\lambda })\Vert J_{a,\lambda }^{\prime }(u_{n})\Vert
_{X_{\lambda }^{-1}}\rightarrow 0\quad \text{as}\ n\rightarrow \infty .
\end{equation*}

\section{Proof of Theorem \protect\ref{t1-1}}

Recall that a $C^{1}$-functional $J_{a,\lambda }$ satisfies Cerami condition
at level $c$ ($(C)_{c}$-condition for short) if any sequence $%
\{u_{n}\}\subset X_{\lambda }$ satisfying
\begin{equation*}
J_{a,\lambda }(u_{n})\rightarrow c\text{ and }(1+\Vert u_{n}\Vert _{\lambda
})\Vert J_{a,\lambda }^{\prime }(u_{n})\Vert _{X_{\lambda }^{-1}}\rightarrow
0,
\end{equation*}%
has a convergent subsequence, and such sequence is called a $(C)_{c}$%
-sequence.

\begin{lemma}
\label{lem3} Assume that $N\geq 1,$ $\delta >0$ and $b>-2A_{0}^{-2}\beta
_{N}^{-1}.$ In addition, assume that conditions $(V1)-(V3),\left( F1\right) $
and $\left( F3\right) $ hold. Then $\{u_{n}\}$ is bounded in $X_{\lambda }$
for each $\lambda \geq \Lambda _{0},$ where $\{u_{n}\}$ is a $(C)_{c}$%
-sequence.
\end{lemma}

\noindent Proof. \/ Following the argument of Lemma $\ref{lem15},$ we can
conclude that the $(C)_{c}$-sequence $\{u_{n}\}$ is bounded in $X_{\lambda }$
for each $\lambda \geq \Lambda _{0}.$ \hfill $\square $

\begin{proposition}
\label{m3} Assume that $b>-2A_{0}^{-2}\beta _{N}^{-1}.$ In addition, we
assume that conditions $(V1)-(V3)$ and $(F1)-(F3)$ hold. Then for each $D>0$%
, there exists $\Lambda _{1}:=\Lambda _{1}(D)\geq \Lambda _{0}>\Lambda _{N}$
such that $J_{a,\lambda }$ satisfies the $(C)_{c}$-condition in $X_{\lambda
} $ for all $c<D$ and $\lambda >\Lambda _{1}.$
\end{proposition}

\noindent Proof. \/ Let $\left\{ u_{n}\right\} $ be a $\left( C\right) _{c}$%
-sequence with $c<D.$ By Lemma \ref{lem3}, $\left\{ u_{n}\right\} $ is
bounded in $X_{\lambda }$ and there exists $D_{0}>0$ such that $\left\Vert
u_{n}\right\Vert _{\lambda }\leq D_{0}.$ Then there exist a subsequence $%
\left\{ u_{n}\right\} $ and $u_{0}$ in $X_{\lambda }$ such that%
\begin{eqnarray*}
u_{n} &\rightharpoonup &u_{0}\text{ weakly in }X_{\lambda }, \\
u_{n} &\rightarrow &u_{0}\text{ strongly in }L_{loc}^{r}(\mathbb{R}^{N}),%
\text{ for }2\leq r<2_{\ast }, \\
u_{n} &\rightarrow &u_{0}\text{ a.e. in }\mathbb{R}^{N}.
\end{eqnarray*}%
Moreover, using $(\ref{13})$ and $(\ref{11})$ implies that the imbedding $%
X_{\lambda }\hookrightarrow W^{1,2}(R^{N})$ is continuous, which shows that%
\begin{equation*}
u_{n}\rightharpoonup u_{0}\text{ weakly in }W^{1,2}(\mathbb{R}^{N}).
\end{equation*}%
Similar to the proof of Lemma 4.4 in \cite{GY}, one can easily obtain that
\begin{equation*}
\nabla u_{n}\left( x\right) \rightarrow \nabla u_{0}\left( x\right) \text{
a.e. in }\mathbb{R}^{N}.
\end{equation*}%
Thus, it follows from Brezis-Lieb lemma \cite{BL} that
\begin{equation}
\int_{\mathbb{R}^{N}}\left\vert \nabla (u_{n}-u_{0})\right\vert ^{2}dx=\int_{%
\mathbb{R}^{N}}\left\vert \nabla u_{n}\right\vert ^{2}dx-\int_{\mathbb{R}%
^{N}}\left\vert \nabla u_{0}\right\vert ^{2}dx+o(1).  \label{18}
\end{equation}%
Now we prove that $u_{n}\rightarrow u_{0}$ strongly in $X_{\lambda }.$ Let $%
v_{n}=u_{n}-u_{0}.$ Then $v_{n}\rightharpoonup 0$ in $X_{\lambda }.$ By the
condition $(V2)$, we have
\begin{equation}
\int_{\mathbb{R}^{N}}v_{n}^{2}dx=\int_{\left\{ V\geq c_{0}\right\}
}v_{n}^{2}dx+\int_{\left\{ V<c_{0}\right\} }v_{n}^{2}dx\leq \frac{1}{\lambda
c_{0}}\left\Vert v_{n}\right\Vert _{\lambda }^{2}+o(1).  \label{20}
\end{equation}%
Using $\left( \ref{20}\right) ,$ together with the H\"{o}lder and Sobolev
inequalities, for any $\lambda >\Lambda _{N},$ we check the following
estimation:

Case $\left( i\right) $ $N=3:$%
\begin{eqnarray}
\int_{\mathbb{R}^{N}}\left\vert v_{n}\right\vert ^{r}dx &\leq &\left\Vert
v_{n}\right\Vert _{L^{\infty }}^{r-2}\int_{\mathbb{R}^{N}}v_{n}^{2}dx\leq
\frac{S_{\infty }^{r-2}}{\lambda c_{0}}\left\Vert v_{n}\right\Vert
_{H^{2}}^{r-2}\left\Vert v_{n}\right\Vert _{\lambda }^{2}+o\left( 1\right)
\notag \\
&\leq &\frac{S_{\infty }^{r-2}}{\lambda c_{0}}\left[ \left( 1+\frac{A_{0}^{2}%
}{2}\right) ^{-1}-S_{\infty }^{2}\left\vert \left\{ V<c_{0}\right\}
\right\vert \right] ^{-(r-2)/2}\left\Vert v_{n}\right\Vert _{\lambda
}^{r}+o(1).  \label{21}
\end{eqnarray}%
Case $\left( ii\right) $\ $N=4:$%
\begin{eqnarray}
\int_{\mathbb{R}^{N}}\left\vert v_{n}\right\vert ^{r}dx &\leq &\left( \int_{%
\mathbb{R}^{N}}v_{n}^{2}dx\right) ^{1/2}\left( \int_{\mathbb{R}%
^{N}}v_{n}^{2\left( r-1\right) }dx\right) ^{1/2}  \notag \\
&\leq &\left( \frac{1}{\lambda c_{0}}\left\Vert v_{n}\right\Vert _{\lambda
}^{2}+o\left( 1\right) \right) ^{1/2}S_{2\left( r-1\right) }^{-\left(
r-1\right) }\left( 1+\frac{A_{0}^{2}}{2}\right) ^{(r-1)/2}\left\Vert
v_{n}\right\Vert _{\lambda }^{r-1}  \notag \\
&=&\frac{S_{2\left( r-1\right) }^{-\left( r-1\right) }}{\sqrt{\lambda c_{0}}}%
\left( 1+\frac{A_{0}^{2}}{2}\right) ^{(r-1)/2}\left\Vert v_{n}\right\Vert
_{\lambda }^{r}+o(1).  \label{22}
\end{eqnarray}%
Case $\left( iii\right) $ $N>4:$%
\begin{eqnarray}
\int_{\mathbb{R}^{N}}\left\vert v_{n}\right\vert ^{r}dx &\leq &\left( \int_{%
\mathbb{R}^{N}}\left\vert v_{n}\right\vert ^{2}dx\right) ^{\frac{2_{\ast }-r%
}{2_{\ast }-2}}\left( \int_{\mathbb{R}^{N}}\left\vert v_{n}\right\vert
^{2_{\ast }}dx\right) ^{\frac{r-2}{2_{\ast }-2}}  \notag \\
&\leq &C_{0}^{\frac{2_{\ast }\left( r-2\right) }{2_{\ast }-2}}\left( \frac{1%
}{\lambda c_{0}}\right) ^{\frac{2_{\ast }-r}{2_{\ast }-2}}\left\Vert
v_{n}\right\Vert _{\lambda }^{r}+o(1).  \label{23}
\end{eqnarray}%
Set
\begin{equation*}
\Psi _{r}:=\left\{
\begin{array}{ll}
\frac{S_{\infty }^{r-2}}{\lambda c_{0}}\left[ \left( 1+\frac{A_{0}^{2}}{2}%
\right) ^{-1}-S_{\infty }^{2}\left\vert \left\{ V<c_{0}\right\} \right\vert %
\right] ^{-(r-1)/2} & \text{ if }N=3, \\
\frac{S_{2\left( r-1\right) }^{-\left( r-1\right) }}{\sqrt{\lambda c_{0}}}%
\left( 1+\frac{A_{0}^{2}}{2}\right) ^{(r-1)/2} & \text{ if }N=4, \\
C_{0}^{\frac{2_{\ast }\left( r-2\right) }{2_{\ast }-2}}\left( \frac{1}{%
\lambda c_{0}}\right) ^{(2_{\ast }-r)/(2_{\ast }-2)} & \text{ if }N>4.%
\end{array}%
\right.
\end{equation*}%
Clearly, $\Psi _{r}\rightarrow 0$ as $\lambda \rightarrow \infty .$ The
inqualities $\left( \ref{21}\right) -\left( \ref{23}\right) $ indicate that%
\begin{equation}
\int_{\mathbb{R}^{N}}\left\vert v_{n}\right\vert ^{r}dx\leq \Psi
_{r}\left\Vert v_{n}\right\Vert _{\lambda }^{r}+o(1).  \label{24}
\end{equation}%
Following the argument of \cite{SW}, it is easy to verify that%
\begin{equation}
\int_{\mathbb{R}^{N}}F(x,v_{n})dx=\int_{\mathbb{R}^{N}}F(x,u_{n})dx-\int_{%
\mathbb{R}^{N}}F(x,u_{0})dx+o(1)  \label{14}
\end{equation}%
and%
\begin{equation*}
\sup_{\left\Vert h\right\Vert _{\lambda }=1}\int_{\mathbb{R}%
^{N}}[f(x,v_{n})-f(x,u_{n})+f(x,u_{0})]h(x)dx=o(1).
\end{equation*}%
Thus, using $(\ref{18}),$ $(\ref{14})$ and Brezis-Lieb Lemma \cite{BL}, we
deduce that%
\begin{eqnarray}
J_{a,\lambda }\left( u_{n}\right) -J_{a,\lambda }\left( u_{0}\right)  &=&%
\frac{1}{2}\left\Vert v_{n}\right\Vert _{\lambda }^{2}+\frac{a}{2\left(
1+\delta \right) }\left( \left\Vert \nabla u_{n}\right\Vert
_{L^{2}}^{2(1+\delta )}-\left\Vert \nabla u_{0}\right\Vert
_{L^{2}}^{2(1+\delta )}\right)   \notag \\
&&+\frac{b}{2}\left\Vert \nabla v_{n}\right\Vert _{L^{2}}^{2}-\int_{\mathbb{R%
}^{N}}F(x,v_{n})dx+o\left( 1\right) .  \label{3.9}
\end{eqnarray}%
Moreover, it follows from the boundedness of the sequence $\left\{
u_{n}\right\} $ in $X_{\lambda }$ and $(\ref{13})$ that there exists a
constant $A>0$ such that%
\begin{equation*}
\left\Vert \nabla u_{n}\right\Vert _{L^{2}}^{2}\rightarrow A\text{ as }%
n\rightarrow \infty .
\end{equation*}%
It indicates that for any $\varphi \in C_{0}^{\infty }(R^{N}),$ there holds%
\begin{eqnarray*}
o(1) &=&\left\langle J_{a,\lambda }^{\prime }\left( u_{n}\right) ,\varphi
\right\rangle  \\
&=&\int_{\mathbb{R}^{N}}\Delta u_{n}\Delta \varphi dx+\int_{\mathbb{R}%
^{N}}\lambda V(x)u_{n}\varphi dx+a\left\Vert \nabla u_{n}\right\Vert
_{L^{2}}^{2\delta }\int_{\mathbb{R}^{N}}\nabla u_{n}\nabla \varphi dx \\
&&+b\int_{\mathbb{R}^{N}}\nabla u_{n}\nabla \varphi dx-\int_{\mathbb{R}%
^{N}}f(x,u_{n})\varphi dx \\
&\rightarrow &\int_{\mathbb{R}^{N}}\Delta u_{0}\Delta \varphi dx+\int_{%
\mathbb{R}^{N}}\lambda V(x)u_{0}\varphi dx+aA^{\delta }\int_{\mathbb{R}%
^{N}}\nabla u_{0}\nabla \varphi dx \\
&&+b\int_{\mathbb{R}^{N}}\nabla u_{0}\nabla \varphi dx-\int_{\mathbb{R}%
^{N}}f(x,u_{0})\varphi dx\text{ as }n\rightarrow \infty ,
\end{eqnarray*}%
which shows that%
\begin{equation*}
\left\Vert u_{0}\right\Vert _{\lambda }^{2}+\left( aA^{\delta }+b\right)
\int_{\mathbb{R}^{N}}\left\vert \nabla u_{0}\right\vert ^{2}dx-\int_{\mathbb{%
R}^{N}}f(x,u_{0})u_{0}dx=o\left( 1\right) .
\end{equation*}%
Note that%
\begin{eqnarray*}
o(1) &=&\left\langle J_{a,\lambda }^{\prime }\left( u_{n}\right)
,u_{n}\right\rangle  \\
&=&\left\Vert u_{n}\right\Vert _{\lambda }^{2}+a\left\Vert \nabla
u_{n}\right\Vert _{L^{2}}^{2(1+\delta )}+b\left\Vert \nabla u_{n}\right\Vert
_{L^{2}}^{2}-\int_{\mathbb{R}^{N}}f(x,u_{n})u_{n}dx.
\end{eqnarray*}%
Combining the above two equalities gives%
\begin{eqnarray}
o\left( 1\right)  &=&\left\Vert u_{n}\right\Vert _{\lambda }^{2}+a\left\Vert
\nabla u_{n}\right\Vert _{L^{2}}^{2(1+\delta )}+b\left\Vert \nabla
u_{n}\right\Vert _{L^{2}}^{2}-\int_{\mathbb{R}^{N}}f(x,u_{n})u_{n}dx  \notag
\\
&&-\left\Vert u_{0}\right\Vert _{\lambda }^{2}-\left( aA^{\delta }+b\right)
\int_{\mathbb{R}^{N}}\left\vert \nabla u_{0}\right\vert ^{2}dx+\int_{\mathbb{%
R}^{N}}f(x,u_{0})u_{0}dx  \notag \\
&=&\left\Vert v_{n}\right\Vert _{\lambda }^{2}+a\left\Vert \nabla
u_{n}\right\Vert _{L^{2}}^{2(1+\delta )}-a\left\Vert \nabla u_{n}\right\Vert
_{L^{2}}^{2\delta }\left\Vert \nabla u_{0}\right\Vert _{L^{2}}^{2}  \notag \\
&&+b\left\Vert \nabla v_{n}\right\Vert _{L^{2}}^{2}-\int_{\mathbb{R}%
^{N}}f(x,v_{n})v_{n}dx+o(1)  \notag \\
&=&\left\Vert v_{n}\right\Vert _{\lambda }^{2}+a\left\Vert \nabla
u_{n}\right\Vert _{L^{2}}^{2\delta }\left\Vert \nabla v_{n}\right\Vert
_{L^{2}}^{2}+b\left\Vert \nabla v_{n}\right\Vert _{L^{2}}^{2}-\int_{\mathbb{R%
}^{N}}f(x,v_{n})v_{n}dx+o(1).  \label{3.10}
\end{eqnarray}%
By Lemma \ref{lem15}, there exists a constant $K<0$ such that%
\begin{equation}
J_{a,\lambda }\left( u_{0}\right) \geq K.  \label{3.11}
\end{equation}%
Thus, in virtue of condition $(F2)$ and $(\ref{3.9})-(\ref{3.11})$ one has%
\begin{eqnarray*}
D-K &\geq &c-J_{a,\lambda }\left( u_{0}\right)  \\
&\geq &J_{a,\lambda }\left( u_{n}\right) -J_{a,\lambda }\left( u_{0}\right)
+o\left( 1\right)  \\
&\geq &\frac{1}{2}\left\Vert v_{n}\right\Vert _{\lambda }^{2}+\frac{a}{%
2\left( 1+\delta \right) }\left( \left\Vert \nabla u_{n}\right\Vert
_{L^{2}}^{2(1+\delta )}-\left\Vert \nabla u_{0}\right\Vert
_{L^{2}}^{2(1+\delta )}\right)  \\
&&+\frac{b}{2}\left\Vert \nabla v_{n}\right\Vert _{L^{2}}^{2}-\int_{\mathbb{R%
}^{N}}F(x,v_{n})dx+o\left( 1\right)  \\
&\geq &\frac{\delta }{2\left( 1+\delta \right) }\left\Vert v_{n}\right\Vert
_{\lambda }^{2}+\frac{b\delta }{2\left( 1+\delta \right) }\left\Vert \nabla
v_{n}\right\Vert _{L^{2}}^{2}-\frac{d_{0}}{2\left( 1+\delta \right) }\int_{%
\mathbb{R}^{N}}v_{n}^{2}dx \\
&&+\frac{a}{2\left( 1+\delta \right) }\left( \left\Vert \nabla
u_{n}\right\Vert _{L^{2}}^{2\delta }-\left\Vert \nabla u_{0}\right\Vert
_{L^{2}}^{2\delta }\right) \left\Vert \nabla u_{0}\right\Vert
_{L^{2}}^{2}+o(1) \\
&\geq &\frac{\delta -d_{0}\Theta _{2,N}^{2}}{2\left( 1+\delta \right) }%
\left\Vert v_{n}\right\Vert _{\lambda }^{2}+\frac{b\delta }{2\left( 1+\delta
\right) }\left\Vert \nabla v_{n}\right\Vert _{L^{2}}^{2} \\
&&+\frac{a}{2\left( 1+\delta \right) }\left( \left\Vert \nabla
u_{n}\right\Vert _{L^{2}}^{2\delta }-\left\Vert \nabla u_{0}\right\Vert
_{L^{2}}^{2\delta }\right) \left\Vert \nabla u_{0}\right\Vert
_{L^{2}}^{2}+o(1),
\end{eqnarray*}%
which implies that there exists a constant $\widehat{D}=\widehat{D}(a,D)>0$
such that%
\begin{equation}
\left\Vert v_{n}\right\Vert _{\lambda }^{2}\leq \widehat{D}+o(1)\text{ for
every }\lambda >\Lambda _{N}.  \label{26}
\end{equation}%
It follows from the condition $(F3),(\ref{24})$ and $(\ref{26})$ that%
\begin{eqnarray*}
o\left( 1\right)  &=&\left\Vert v_{n}\right\Vert _{\lambda }^{2}+a\left\Vert
\nabla u_{n}\right\Vert _{L^{2}}^{2\delta }\left\Vert \nabla
v_{n}\right\Vert _{L^{2}}^{2}+b\left\Vert \nabla v_{n}\right\Vert
_{L^{2}}^{2}-\int_{\mathbb{R}^{N}}f(x,v_{n})v_{n}dx \\
&\geq &\left\Vert v_{n}\right\Vert _{\lambda }^{2}+b\left\Vert \nabla
v_{n}\right\Vert _{L^{2}}^{2}-\epsilon \int_{\mathbb{R}%
^{N}}v_{n}^{2}dx-C_{1,\epsilon }\int_{\mathbb{R}^{N}}\left\vert
v_{n}\right\vert ^{p}dx \\
&\geq &\left\{
\begin{array}{ll}
\Vert v_{n}\Vert _{\lambda }^{2}-\epsilon \Psi _{r}^{2}\widehat{D}%
-C_{1,\epsilon }\Psi _{r}^{p}\widehat{D}^{p/2} & \text{ if }b\geq 0, \\
\frac{1}{2}\left( 2+bA_{0}^{2}\beta _{N}\right) \Vert v_{n}\Vert _{\lambda
}^{2}-\epsilon \Psi _{r}^{2}\widehat{D}-C_{1,\epsilon }\Psi _{r}^{p}\widehat{%
D}^{p/2} & \text{ if }-2A_{0}^{-2}\beta _{N}^{-1}<b<0,%
\end{array}%
\right.
\end{eqnarray*}%
which implies that there exists $\Lambda _{1}:=\Lambda _{1}(a,D)\geq \Lambda
_{0}>\Lambda _{N}$ such that for each $\lambda >\Lambda _{1},$%
\begin{equation*}
v_{n}\rightarrow 0\text{ strongly in }X_{\lambda }.
\end{equation*}%
This completes the proof. \hfill $\square $

\begin{theorem}
\label{t3}Assume that $N\geq 3,$ $\delta \geq \frac{2}{N-2}$ and $%
b>-2A_{0}^{-2}\beta _{N}^{-1}.$ In addition, we assume that conditions $%
(V1)-(V3)$ and $(F1)-(F3)$ are satisfied. Then for each $0<a<a_{\ast }$ and $%
\lambda >\Lambda _{1},$ $J_{a,\lambda }$ has a nonzero critical point $%
u_{a,\lambda }^{(1)}\in X_{\lambda }$ such that $J_{a,\lambda }\left(
u_{a,\lambda }^{(1)}\right) =c_{\lambda }>0.$
\end{theorem}

\noindent Proof. \/ By virtue of Theorem \ref{l2}, Lemmas \ref{lem1} and \ref%
{lem2}, for every $\lambda >\Lambda _{1}$ and $0<a<a^{\ast },$ there exists
a sequence $\{u_{n}\}\subset X_{\lambda }$ satisfying%
\begin{equation*}
J_{a,\lambda }(u_{n})\rightarrow c_{\lambda }>0\quad \text{and}\quad
(1+\Vert u_{n}\Vert _{\lambda })\Vert J_{a,\lambda }^{\prime }(u_{n})\Vert
_{X_{\lambda }^{-1}}\rightarrow 0,\quad \text{as}\ n\rightarrow \infty .
\end{equation*}%
By Lemma \ref{lem3}, one has $\{u_{n}\}$ is bounded in $X_{\lambda }.$ Then
it follows from Proposition \ref{m3} and the fact of $0<\eta \leq c_{\lambda
}\leq c_{0}\left( \Omega \right) $ that $J_{a,\lambda }$ satisfies the (C)$%
_{\alpha }$--condition in $X_{_{\lambda }}$ for all $c_{\lambda }<D$ and $%
\lambda >\Lambda _{1}.$ This indicates that there exist a subsequence $%
\{u_{n}\}$ and $u_{a,\lambda }^{(1)}\in X_{\lambda }$ such that $%
u_{n}\rightarrow u_{a,\lambda }^{(1)}$ strongly in $X_{\lambda }.$ The proof
is completed. \hfill $\square $

\begin{lemma}
\label{lem5}Suppose that $N\geq 3,$ $\delta \geq \frac{2}{N-2}$ and $%
b>-2A_{0}^{-2}\beta _{N}^{-1}.$ In addition, assume that conditions $%
(V1)-(V3),(F1)$ and $(F3)$ hold. Then for every $0<a<a_{\ast }$ and $\lambda
>\Lambda _{1}$ one has
\begin{equation}
-\infty <\theta _{a}=:\inf \left\{ J_{a,\lambda }(u):u\in X_{\lambda }\text{
with }\rho <\left\Vert u\right\Vert _{\lambda }<R_{a}\right\} <\frac{\kappa
}{2}<0.  \label{3.8}
\end{equation}
\end{lemma}

\noindent Proof. \/ The proof directly follows from Lemmas \ref{lem15} and %
\ref{lem2}. \hfill $\square $

\begin{theorem}
\label{t5}Suppose that $N\geq 3,$ $\delta \geq \frac{2}{N-2}$ and $%
b>-2A_{0}^{-2}\beta _{N}^{-1}.$ In addition, assume that conditions $%
(V1)-(V3),(F1)-(F3)$ hold. Then for every $0<a<a_{\ast }$ and $\lambda
>\Lambda _{1},$ $J_{a,\lambda }$ has a nonzero critical point $u_{a,\lambda
}^{(2)}\in X_{\lambda }$ such that
\begin{equation*}
J_{a,\lambda }\left( u_{a,\lambda }^{(2)}\right) =\theta _{a}<0,
\end{equation*}%
where $\widehat{\theta }_{a}$ is as in $(\ref{3.8})$. Furthermore, when $%
\delta >\frac{2}{N-2},$ for every $\lambda >\Lambda _{1}$ there holds%
\begin{equation*}
J_{a,\lambda }\left( u_{a,\lambda }^{(2)}\right) \rightarrow -\infty \ \text{%
and}\ \left\Vert u_{a,\lambda }^{(2)}\right\Vert _{\lambda }\rightarrow
\infty \text{ as }a\rightarrow 0,
\end{equation*}
\end{theorem}

\noindent Proof. \/ It follows from Lemmas \ref{lem3}, \ref{lem5} and the
Ekeland variational principle that there exists a minimizing bounded
sequence $\{u_{n}\}\subset X_{\lambda }$ with $\rho <\Vert u_{n}\Vert
_{\lambda }<R_{a}$ such that%
\begin{equation*}
J_{a,\lambda }(u_{n})\rightarrow \theta _{a}\text{ and }J_{a,\lambda
}^{\prime }(u_{n})\rightarrow 0\text{ as }n\rightarrow \infty .
\end{equation*}%
Similar to the proof of Theorem \ref{t3}, there exist a subsequence $%
\{u_{n}\}$ and $u_{a,\lambda }^{(2)}\in X_{\lambda }$ with $\rho <\left\Vert
u_{a,\lambda }^{(2)}\right\Vert _{\lambda }<R_{a}$ such that $%
u_{n}\rightarrow u_{a,\lambda }^{(2)}$ strongly in $X_{\lambda },$ which
implies that $J_{a,\lambda }^{\prime }\left( u_{a,\lambda }^{(2)}\right) =0$
and $J_{a,\lambda }\left( u_{a,\lambda }^{(2)}\right) =\theta _{a}<0.$
Furthermore, by Lemmas \ref{lem15} and \ref{lem2} we have
\begin{equation*}
J_{a,\lambda }\left( u_{a,\lambda }^{(2)}\right) \leq J_{a,\lambda
}(e)\rightarrow -\infty \ \text{as}\ a\rightarrow 0.
\end{equation*}%
It implies that
\begin{equation*}
\left\Vert u_{a,\lambda }^{(2)}\right\Vert _{\lambda }\rightarrow \infty
\text{ as }a\rightarrow 0,
\end{equation*}%
Consequently, we complete the proof. \hfill $\square $

We are now ready to prove Theorem \ref{t1-1}. By Theorems \ref{t3} and \ref%
{t5}, for every $0<a<a_{\ast }$ and $\lambda >\Lambda _{1},$ there exist two
nontrivial solutions $u_{a,\lambda }^{(1)}$ and $u_{a,\lambda }^{(2)}$ of
Eq. $(K_{a,\lambda })$ such that%
\begin{equation*}
J_{a,\lambda }\left( u_{a,\lambda }^{(2)}\right) =\theta _{a}<\frac{\kappa }{%
2}<0<\eta <c_{\lambda }=J_{a,\lambda }\left( u_{a,\lambda }^{(1)}\right) ,
\end{equation*}%
which implies that $u_{a,\lambda }^{(1)}\neq u_{a,\lambda }^{(2)}.$
Furthermore, when $\delta >\frac{2}{N-2},$ for every $\lambda >\Lambda _{1}$
there holds%
\begin{equation*}
J_{a,\lambda }\left( u_{a,\lambda }^{(2)}\right) \rightarrow -\infty \ \text{%
and}\ \left\Vert u_{a,\lambda }^{(2)}\right\Vert _{\lambda }\rightarrow
\infty \text{ as }a\rightarrow 0.
\end{equation*}%
Since $J_{a,\lambda }(u)\geq 0$ on $\left\{ u\in X_{\lambda }\text{ with }%
\left\Vert u\right\Vert _{\lambda }\leq \rho \cup \left\Vert u\right\Vert
_{\lambda }\geq R_{a}\right\} $ by Lemmas \ref{lem5} and \ref{lem1}, we
conclude that $u_{a,\lambda }^{(2)}$ is a ground state solution of Eq. $%
(K_{a,\lambda }).$ This completes the proof of Theorem \ref{t1-1}.

\section{Proof of Theorem \protect\ref{t1-2}}

In this section, we give the proof of Theorem \ref{t1-2}: Let $u_{0}$ be a
nontrivial solution of Eq. $(K_{a,\lambda }).$ Then there holds
\begin{equation*}
\left\Vert u_{0}\right\Vert _{\lambda }^{2}+a\left\Vert \nabla
u_{0}\right\Vert _{L^{2}}^{2(1+\delta )}+b\left\Vert \nabla u_{0}\right\Vert
_{L^{2}}^{2}-\int_{\mathbb{R}^{N}}f\left( x,u_{0}\right) u_{0}dx=0.
\end{equation*}%
We now divide the proof into two separate cases:

Case $A:\int_{\mathbb{R}^{N}}\lambda V(x)u_{0}^{2}dx\geq \lambda c_{0}\left(
\frac{C_{1,\epsilon }}{\lambda c_{0}-\epsilon }\right) ^{\frac{4}{(p-2)(N-2)}%
}\left( \overline{S}^{-1}\left\Vert \nabla u_{0}\right\Vert _{L^{2}}\right)
^{\frac{2N}{N-2}}.$ It follows from the condition $(F3)$ and $(\ref{8})$ that%
\begin{eqnarray*}
0 &=&\left\Vert u_{0}\right\Vert _{\lambda }^{2}+a\left\Vert \nabla
u_{0}\right\Vert _{L^{2}}^{2(1+\delta )}+b\left\Vert \nabla u_{0}\right\Vert
_{L^{2}}^{2}-\int_{\mathbb{R}^{N}}f\left( x,u_{0}\right) u_{0}dx \\
&\geq &\left\Vert u_{0}\right\Vert _{\lambda }^{2}+a\left\Vert \nabla
u_{0}\right\Vert _{L^{2}}^{2(1+\delta )}+b\left\Vert \nabla u_{0}\right\Vert
_{L^{2}}^{2} \\
&&-\epsilon \left( \frac{1}{\lambda c_{0}}\int_{\mathbb{R}^{N}}\lambda
V\left( x\right) u_{0}^{2}dx+\overline{S}^{-2}\left\vert \left\{
V<c_{0}\right\} \right\vert ^{\frac{2}{N}}\left\Vert \nabla u_{0}\right\Vert
_{L^{2}}^{2}\right) \\
&&-\frac{C_{1,\epsilon }}{\overline{S}^{\frac{N\left( p-2\right) }{2}}}%
\left( \frac{1}{\lambda c_{0}}\int_{\mathbb{R}^{N}}\lambda V\left( x\right)
u_{0}^{2}dx+\overline{S}^{-2}\left\vert \left\{ V<c_{0}\right\} \right\vert
^{\frac{2}{N}}\left\Vert \nabla u_{0}\right\Vert _{L^{2}}^{2}\right) ^{\frac{%
2N-p(N-2)}{4}}\left\Vert \nabla u_{0}\right\Vert _{L^{2}}^{\frac{N\left(
p-2\right) }{2}} \\
&\geq &a\left\Vert \nabla u_{0}\right\Vert _{L^{2}}^{2(1+\delta )}+\left(
b-\epsilon \overline{S}^{-2}\left\vert \left\{ V<c_{0}\right\} \right\vert ^{%
\frac{2}{N}}\right) \left\Vert \nabla u_{0}\right\Vert _{L^{2}}^{2} \\
&&-\frac{C_{1,\epsilon }}{\overline{S}^{p}}\left\vert \left\{
V<c_{0}\right\} \right\vert ^{\frac{2N-p\left( N-2\right) }{2N}}\left\Vert
\nabla u_{0}\right\Vert _{L^{2}}^{p}+\left( \int_{\mathbb{R}^{N}}\lambda
V(x)u_{0}^{2}dx\right) ^{\frac{2N-p\left( N-2\right) }{4}} \\
&&\cdot \left[ \frac{\lambda c_{0}-\epsilon }{\lambda c_{0}}\left( \int_{%
\mathbb{R}^{N}}\lambda V(x)u^{2}dx\right) ^{\frac{(p-2)(N-2)}{4}}-\frac{%
C_{1,\epsilon }}{\overline{S}^{\frac{N\left( p-2\right) }{2}}}\left( \frac{1%
}{\lambda c_{0}}\right) ^{\frac{2N-p\left( N-2\right) }{4}}\left\Vert \nabla
u_{0}\right\Vert _{L^{2}}^{\frac{N\left( p-2\right) }{2}}\right] \\
&\geq &a\left\Vert \nabla u_{0}\right\Vert _{L^{2}}^{2(1+\delta )}+\left(
b-\epsilon \overline{S}^{-2}\left\vert \left\{ V<c_{0}\right\} \right\vert ^{%
\frac{2}{N}}\right) \left\Vert \nabla u_{0}\right\Vert _{L^{2}}^{2}-\frac{%
C_{1,\epsilon }}{\overline{S}^{p}}\left\vert \left\{ V<c_{0}\right\}
\right\vert ^{\frac{2N-p\left( N-2\right) }{2N}}\left\Vert \nabla
u_{0}\right\Vert _{L^{2}}^{p} \\
&>&0.
\end{eqnarray*}%
provided that%
\begin{equation*}
a>\frac{p-2}{2\delta }\left[ \frac{2(\delta +1)-p}{2\delta \left( b-\epsilon
\overline{S}^{-2}\left\vert \left\{ V<c_{0}\right\} \right\vert
^{2/N}\right) }\right] ^{\frac{2(\delta +1)-p}{p-2}}\left( \frac{%
C_{1,\epsilon }}{\overline{S}^{p}}\left\vert \left\{ V<c_{0}\right\}
\right\vert ^{\frac{2N-p\left( N-2\right) }{2N}}\right) ^{\frac{2\delta }{p-2%
}}.
\end{equation*}%
This is a contradiction.

Case $B:\int_{\mathbb{R}^{N}}\lambda V(x)u_{0}^{2}dx<\lambda c_{0}\left[
\frac{C_{1,\epsilon }}{(\lambda c_{0}-\epsilon )}\right] ^{\frac{4}{%
(p-2)(N-2)}}\left( \overline{S}^{-1}\left\Vert \nabla u_{0}\right\Vert
_{L^{2}}\right) ^{\frac{2N}{N-2}}.$ By virtue of $(\ref{8})$ one has%
\begin{eqnarray*}
&&\int_{\mathbb{R}^{N}}|u_{0}|^{p}dx \\
&\leq &\left( \frac{1}{\lambda c_{0}}\int_{\mathbb{R}^{N}}\lambda
V(x)u_{0}^{2}dx+\frac{\left\vert \left\{ V<c_{0}\right\} \right\vert ^{\frac{%
2}{N}}}{\overline{S}^{2}}\left\Vert \nabla u_{0}\right\Vert
_{L^{2}}^{2}\right) ^{\frac{2^{\ast }-p}{2^{\ast }-2}}\left( \overline{S}%
^{-1}\left\Vert \nabla u_{0}\right\Vert _{L^{2}}\right) ^{\frac{N\left(
p-2\right) }{2}} \\
&<&\overline{S}^{-2^{\ast }}\left( \frac{C_{1,\epsilon }}{\lambda
c_{0}-\epsilon }\right) ^{\frac{2N-p\left( N-2\right) }{(p-2)(N-2)}%
}\left\Vert \nabla u_{0}\right\Vert _{L^{2}}^{2^{\ast }}+\overline{S}%
^{-p}\left\vert \left\{ V<c_{0}\right\} \right\vert ^{\frac{2N-p\left(
N-2\right) }{2N}}\left\Vert \nabla u_{0}\right\Vert _{L^{2}}^{p}.
\end{eqnarray*}%
Using this, together with condition $(F3),$ gives%
\begin{eqnarray*}
0 &=&\left\Vert u_{0}\right\Vert _{\lambda }^{2}+a\left\Vert \nabla
u_{0}\right\Vert _{L^{2}}^{2(1+\delta )}+b\left\Vert \nabla u_{0}\right\Vert
_{L^{2}}^{2}-\int_{\mathbb{R}^{N}}f\left( x,u_{0}\right) u_{0}dx \\
&>&a\left\Vert \nabla u_{0}\right\Vert _{L^{2}}^{2(1+\delta )}+\left(
b-\epsilon \overline{S}^{-2}\left\vert \left\{ V<c_{0}\right\} \right\vert ^{%
\frac{2}{N}}\right) \left\Vert \nabla u_{0}\right\Vert _{L^{2}}^{2} \\
&&-\frac{C_{1,\epsilon }^{\frac{4}{(p-2)(N-2)}}}{\overline{S}^{2^{\ast }}}%
(\lambda c_{0}-\epsilon )^{-\frac{2N-p\left( N-2\right) }{(p-2)(N-2)}%
}\left\Vert \nabla u_{0}\right\Vert _{L^{2}}^{2^{\ast }}-\frac{C_{1,\epsilon
}}{\overline{S}^{p}}\left\vert \left\{ V<c_{0}\right\} \right\vert ^{\frac{%
2N-p\left( N-2\right) }{2N}}\left\Vert \nabla u_{0}\right\Vert _{L^{2}}^{p}.
\end{eqnarray*}%
If $\delta =\frac{2}{N-2},$ then for%
\begin{eqnarray*}
a &>&\frac{p-2}{2\delta }\left[ \frac{2(\delta +1)-p}{2\delta \left(
b-\epsilon \overline{S}^{-2}\left\vert \left\{ V<c_{0}\right\} \right\vert
^{2/N}\right) }\right] ^{\frac{2(\delta +1)-p}{p-2}}\left( \frac{%
C_{1,\epsilon }}{\overline{S}^{p}}\left\vert \left\{ V<c_{0}\right\}
\right\vert ^{\frac{2N-p\left( N-2\right) }{2N}}\right) ^{\frac{2\delta }{p-2%
}} \\
&&+\frac{C_{1,\epsilon }}{\overline{S}^{2^{\ast }}}\left( \frac{%
C_{1,\epsilon }}{\lambda c_{0}-\epsilon }\right) ^{\frac{2N-p\left(
N-2\right) }{(p-2)(N-2)}},
\end{eqnarray*}%
there holds%
\begin{eqnarray*}
0 &>&a\left\Vert \nabla u_{0}\right\Vert _{L^{2}}^{2(1+\delta )}+\left(
b-\epsilon \overline{S}^{-2}\left\vert \left\{ V<c_{0}\right\} \right\vert ^{%
\frac{2}{N}}\right) \left\Vert \nabla u_{0}\right\Vert _{L^{2}}^{2} \\
&&-\frac{C_{1,\epsilon }}{\overline{S}^{2^{\ast }}}\left( \frac{%
C_{1,\epsilon }}{\lambda c_{0}-\epsilon }\right) ^{\frac{2N-p\left(
N-2\right) }{(p-2)(N-2)}}\left\Vert \nabla u_{0}\right\Vert
_{L^{2}}^{2^{\ast }}-\frac{C_{1,\epsilon }}{\overline{S}^{p}}\left\vert
\left\{ V<c_{0}\right\} \right\vert ^{\frac{2N-p\left( N-2\right) }{2N}%
}\left\Vert \nabla u_{0}\right\Vert _{L^{2}}^{p} \\
&>&0.
\end{eqnarray*}%
This is a contradiction. If $\delta >\frac{2}{N-2},$ then we consider the
following two cases:\newline
$(i)$ $\left\Vert \nabla u_{0}\right\Vert _{L^{2}}^{2}\geq \overline{S}%
^{2}\left\vert \left\{ V<c_{0}\right\} \right\vert ^{\frac{N-2}{N}}\left(
\frac{\lambda c_{0}-\epsilon }{C_{1,\epsilon }}\right) ^{\frac{2}{p-2}}.$
Then we have%
\begin{eqnarray*}
0 &>&a\left\Vert \nabla u_{0}\right\Vert _{L^{2}}^{2(1+\delta )}+\left(
b-\epsilon \overline{S}^{-2}\left\vert \left\{ V<c_{0}\right\} \right\vert ^{%
\frac{2}{N}}\right) \left\Vert \nabla u_{0}\right\Vert _{L^{2}}^{2} \\
&&-\frac{C_{1,\epsilon }}{\overline{S}^{2^{\ast }}}\left( \frac{%
C_{1,\epsilon }}{(\lambda c_{0}-\epsilon )}\right) ^{\frac{2N-p\left(
N-2\right) }{(p-2)(N-2)}}\left\Vert \nabla u_{0}\right\Vert
_{L^{2}}^{2^{\ast }}-\frac{C_{1,\epsilon }}{\overline{S}^{p}}\left\vert
\left\{ V<c_{0}\right\} \right\vert ^{\frac{2N-p\left( N-2\right) }{2N}%
}\left\Vert \nabla u_{0}\right\Vert _{L^{2}}^{p} \\
&\geq &a\left\Vert \nabla u_{0}\right\Vert _{L^{2}}^{2(1+\delta )}+\left(
b-\epsilon \overline{S}^{-2}\left\vert \left\{ V<c_{0}\right\} \right\vert ^{%
\frac{2}{N}}\right) \left\Vert \nabla u_{0}\right\Vert _{L^{2}}^{2} \\
&&-\frac{2C_{1,\epsilon }}{\overline{S}^{2^{\ast }}}\left( \frac{%
C_{1,\epsilon }}{(\lambda c_{0}-\epsilon )}\right) ^{\frac{2N-p\left(
N-2\right) }{(p-2)(N-2)}}\left\Vert \nabla u_{0}\right\Vert
_{L^{2}}^{2^{\ast }} \\
&>&0,
\end{eqnarray*}%
provided that%
\begin{equation*}
a>\frac{2^{\frac{2+\delta (N-2)}{2}}C_{1,\epsilon }^{\frac{2\delta }{p-2}%
}(\lambda c_{0}-\epsilon )^{-\frac{[2N-p\left( N-2\right) ]\delta }{2(p-2)}}%
}{\delta (N-2)\overline{S}^{N\delta }}\left[ \frac{\delta (N-2)-2}{\delta
(N-2)\left( b-\epsilon \overline{S}^{-2}\left\vert \left\{ V<c_{0}\right\}
\right\vert ^{\frac{2}{N}}\right) }\right] ^{\frac{\delta (N-2)-2}{2}}.
\end{equation*}%
This is a contradiction.\newline
$(ii)$ $\left\Vert \nabla u_{0}\right\Vert _{L^{2}}^{2}<\overline{S}%
^{2}\left\vert \left\{ V<c_{0}\right\} \right\vert ^{\frac{N-2}{N}}\left(
\frac{\lambda c_{0}-\epsilon }{C_{1,\epsilon }}\right) ^{\frac{2}{p-2}}.$
There holds%
\begin{eqnarray*}
0 &>&a\left\Vert \nabla u_{0}\right\Vert _{L^{2}}^{2(1+\delta )}+\left(
b-\epsilon \overline{S}^{-2}\left\vert \left\{ V<c_{0}\right\} \right\vert ^{%
\frac{2}{N}}\right) \left\Vert \nabla u_{0}\right\Vert _{L^{2}}^{2} \\
&&-\frac{C_{1,\epsilon }}{\overline{S}^{2^{\ast }}}\left( \frac{%
C_{1,\epsilon }}{(\lambda c_{0}-\epsilon )}\right) ^{\frac{2N-p\left(
N-2\right) }{(p-2)(N-2)}}\left\Vert \nabla u_{0}\right\Vert
_{L^{2}}^{2^{\ast }}-\frac{C_{1,\epsilon }}{\overline{S}^{p}}\left\vert
\left\{ V<c_{0}\right\} \right\vert ^{\frac{2N-p\left( N-2\right) }{2N}%
}\left\Vert \nabla u_{0}\right\Vert _{L^{2}}^{p} \\
&\geq &a\left\Vert \nabla u_{0}\right\Vert _{L^{2}}^{2(1+\delta )}+\left(
b-\epsilon \overline{S}^{-2}\left\vert \left\{ V<c_{0}\right\} \right\vert ^{%
\frac{2}{N}}\right) \left\Vert \nabla u_{0}\right\Vert _{L^{2}}^{2} \\
&&-2\frac{C_{1,\epsilon }}{\overline{S}^{p}}\left\vert \left\{
V<c_{0}\right\} \right\vert ^{\frac{2N-p\left( N-2\right) }{2N}}\left\Vert
\nabla u_{0}\right\Vert _{L^{2}}^{p} \\
&>&0,
\end{eqnarray*}%
provided that%
\begin{equation*}
a>\frac{p-2}{2\delta }\left[ \frac{2(\delta +1)-p}{2\delta \left( b-\epsilon
\overline{S}^{-2}\left\vert \left\{ V<c_{0}\right\} \right\vert
^{2/N}\right) }\right] ^{\frac{2(\delta +1)-p}{p-2}}\left( \frac{%
2C_{1,\epsilon }\left\vert \left\{ V<c_{0}\right\} \right\vert ^{\frac{%
2N-p\left( N-2\right) }{2N}}}{\overline{S}^{p}}\right) ^{\frac{2\delta }{p-2}%
}.
\end{equation*}%
We also get a contradiction. Therefore, there exists a constant $a^{\ast }>0$
such that for every $a>a^{\ast }$, Eq. $(K_{a,\lambda })$ does not admit any
nontrivial solution for all $\lambda >\Lambda _{N}.$ This completes the
proof of Theorem \ref{t1-2}.

\section{Concentration of solutions}

In this section, we investigate the concentration for solutions and give the
proof of Theorem \ref{t1-3}.

Proof of Theorem \ref{t1-3}: Following the arguments of \cite{BPW}. We
firstly choose a positive sequence $\{\lambda _{n}\}$ such that $\Lambda
_{1}<\lambda _{1}\leq \lambda _{2}\leq ...\leq \lambda _{n}\rightarrow
\infty $ as $n\rightarrow \infty .$ Let $u_{n}^{(i)}:=u_{a,\lambda
_{n}}^{(i)}$ with $i=1,2$ be the critical points of $J_{a,\lambda _{n}}$
obtained in Theorem \ref{t1-1}. Since%
\begin{equation}
J_{a,\lambda _{n}}\left( u_{n}^{(2)}\right) <\frac{\kappa }{2}<0<\eta
<c_{\lambda _{n}}=J_{a,\lambda _{n}}\left( u_{n}^{(1)}\right) <D,
\label{5-3}
\end{equation}%
by Lemma \ref{lem15} one has
\begin{equation}
\left\Vert u_{n}^{(i)}\right\Vert _{\lambda _{n}}\leq C_{0},  \label{5-1}
\end{equation}%
where the constant $C_{0}>0$ is independent of $\lambda _{n}$. This implies
that $\left\Vert u_{n}^{(i)}\right\Vert _{\lambda _{1}}\leq C_{0}.$ Thus,
there exist $u_{\infty }^{(i)}\in X$ $(i=1,2)$ such that%
\begin{eqnarray*}
u_{n}^{(i)} &\rightharpoonup &u_{\infty }^{(i)}\text{ weakly in }X_{\lambda
_{1}}, \\
u_{n}^{(i)} &\rightarrow &u_{\infty }^{(i)}\text{ strongly in }L_{loc}^{r}(%
\mathbb{R}^{N}),\text{ for }2\leq r<2_{\ast }, \\
u_{n}^{(i)} &\rightarrow &u_{\infty }^{(i)}\text{ a.e. in }\mathbb{R}^{N}.
\end{eqnarray*}%
Following the proof of Proposition \ref{m3}, we can conclude that%
\begin{equation*}
u_{n}^{(i)}\rightarrow u_{\infty }^{(i)}\text{ strongly in }X_{\lambda _{1}}.
\end{equation*}%
This shows that $u_{n}^{(i)}\rightarrow u_{\infty }^{(i)}$ strongly in $%
H^{2}(R^{N})$ by (\ref{12}) and that
\begin{equation}
\left\Vert \nabla u_{n}^{(i)}\right\Vert _{L^{2}}^{2}\rightarrow \left\Vert
\nabla u_{\infty }^{(i)}\right\Vert _{L^{2}}^{2}\text{ as }n\rightarrow
\infty  \label{5-2}
\end{equation}%
by (\ref{13}) and (\ref{18}).

Using Fatou's Lemma leads to
\begin{equation*}
\int_{\mathbb{R}^{N}}V(x)\left( u_{\infty }^{(i)}\right) ^{2}dx\leq
\liminf_{n\rightarrow \infty }\int_{\mathbb{R}^{N}}V(x)\left(
u_{n}^{(i)}\right) ^{2}dx\leq \liminf_{n\rightarrow \infty }\frac{\left\Vert
u_{n}^{(i)}\right\Vert _{\lambda _{n}}^{2}}{\lambda _{n}}=0,
\end{equation*}%
which implies that $u_{\infty }^{(i)}(x)=0$ a.e. in $R^{N}\backslash
\overline{\Omega }.$ Moreover, fixing $\phi \in C_{0}^{\infty
}(R^{N}\backslash \overline{\Omega })$, we have
\begin{equation*}
\int_{\mathbb{R}^{N}\backslash \overline{\Omega }}\nabla u_{\infty
}^{(i)}(x)\phi (x)dx=-\int_{\mathbb{R}^{N}\backslash \overline{\Omega }%
}u_{\infty }^{(i)}(x)\nabla \phi (x)dx=0.
\end{equation*}%
This indicates that
\begin{equation*}
\nabla u_{\infty }^{(i)}(x)=0\ \mbox{a.e. in}\ \mathbb{R}^{N}\backslash
\overline{\Omega }.
\end{equation*}%
Since $\partial \Omega $ is smooth, $u_{\infty }^{(i)}\in
H^{2}(R^{N}\backslash \overline{\Omega })$ and $\nabla u_{\infty }^{(i)}\in
H^{1}(R^{N}\backslash \overline{\Omega })$, it follows from Trace Theorem
that there are constants $\overline{C},\widetilde{C}>0$ such that
\begin{equation*}
\Vert u_{\infty }^{(i)}\Vert _{L^{2}(\partial \Omega )}\leq \overline{C}%
\Vert u_{\infty }^{(i)}\Vert _{H^{2}(\mathbb{R}^{N}\backslash \overline{%
\Omega })}=0.
\end{equation*}%
and
\begin{equation*}
\Vert \nabla u_{\infty }^{(i)}\Vert _{L^{2}(\partial \Omega )}\leq
\widetilde{C}\Vert \nabla u_{\infty }^{(i)}\Vert _{H^{1}(\mathbb{R}%
^{N}\backslash \overline{\Omega })}=0.
\end{equation*}%
These show that $u_{\infty }^{(i)}\in H_{0}^{2}(\Omega )$.

Since $\left\langle J_{a,\lambda _{n}}^{\prime }\left( u_{n}^{(i)}\right)
,\varphi \right\rangle =0$ for any $\varphi \in C_{0}^{\infty }\left( \Omega
\right) $, combining (\ref{5-2}), it is not difficult to check that
\begin{equation*}
\int_{\Omega }\Delta u_{\infty }^{(i)}\Delta \varphi dx+\left[ a\left(
\int_{\Omega }\left\vert \nabla u_{\infty }^{(i)}\right\vert ^{2}dx\right)
^{\delta }+b\right] \int_{\Omega }\nabla u_{\infty }^{(i)}\cdot \nabla
\varphi dx=\int_{\Omega }f\left( x,u_{\infty }^{(i)}\right) \varphi dx,
\end{equation*}%
that is, $u_{\infty }^{(i)}$ are the weak solutions of the equation%
\begin{equation*}
\Delta ^{2}u-M\left( \int_{\Omega }|\nabla u|^{2}dx\right) \Delta u=f\left(
x,u\right) \text{ in }\Omega ,
\end{equation*}%
where $M\left( t\right) =at^{\delta }+b.$ Since $u_{n}^{(i)}\rightarrow
u_{\infty }^{(i)}$ strongly in $X,$ using (\ref{5-3}) and the fact that $%
\eta $ and $\kappa $ are independent of $\lambda _{n}$ gives%
\begin{equation*}
\frac{1}{2}\int_{\Omega }\left\vert \Delta u_{\infty }^{(1)}\right\vert
^{2}dx+\frac{a}{2\left( 1+\delta \right) }\left( \int_{\Omega }\left\vert
\nabla u_{\infty }^{(1)}\right\vert ^{2}dx\right) ^{\delta +1}+\frac{b}{2}%
\int_{\Omega }\left\vert \nabla u_{\infty }^{(1)}\right\vert
^{2}dx-\int_{\Omega }F\left( x,u_{\infty }^{(1)}\right) dx\geq \eta >0
\end{equation*}%
and%
\begin{equation*}
\frac{1}{2}\int_{\Omega }\left\vert \Delta u_{\infty }^{(2)}\right\vert
^{2}dx+\frac{a}{2\left( 1+\delta \right) }\left( \int_{\Omega }\left\vert
\nabla u_{\infty }^{(2)}\right\vert ^{2}dx\right) ^{\delta +1}+\frac{b}{2}%
\int_{\Omega }\left\vert \nabla u_{\infty }^{(2)}\right\vert
^{2}dx-\int_{\Omega }F\left( x,u_{\infty }^{(2)}\right) dx\leq \frac{\kappa
}{2}<0.
\end{equation*}%
These imply that $u_{\infty }^{(i)}\neq 0(i=1,2)$ and $u_{\infty
}^{(1)}\not=u_{\infty }^{(2)}.$ Consequently, this complete the proof.

\section{Appendix A}

Consider the following biharmonic equations

\begin{equation}
\left\{
\begin{array}{ll}
\Delta ^{2}u+\lambda V\left( x\right) u=|u|^{p-2}u & \text{ in }\mathbb{R}%
^{N}, \\
u\in H^{2}(\mathbb{R}^{N}), &
\end{array}%
\right.  \label{6-1}
\end{equation}%
where $N\geq 3,2<p<\frac{2N}{N-2},\lambda >0$ is a parameter and $V(x)$
satisfies conditions $(V1)-(V3).$

Eq. $(\ref{6-1})$ is variational and its solutions are critical points of
the functional defined in $X_{\lambda }$ by
\begin{equation*}
\mathcal{J}_{\lambda }\left( u\right) =\frac{1}{2}\left\Vert u\right\Vert
_{\lambda }^{2}-\frac{1}{p}\int_{\mathbb{R}^{N}}|u|^{p}dx,
\end{equation*}%
where $\left\Vert u\right\Vert _{\lambda }$ is defined as (\ref{10}). It is
easily seen that the functional $J_{\lambda }$ is of class $C^{1}$ in $%
X_{\lambda }$, and that
\begin{equation}
\langle \mathcal{J}_{\lambda }^{\prime }(u),v\rangle =\int_{\mathbb{R}^{N}}%
\left[ \Delta u\cdot \Delta v+\lambda V\left( x\right) uv\right] dx-\int_{%
\mathbb{R}^{N}}|u|^{p-2}uvdx.  \notag
\end{equation}

Define the Nehari manifold by
\begin{equation*}
\mathcal{N}=\left\{ u\in X_{\lambda }\backslash \{0\}\text{ }|\text{ }%
\langle \mathcal{J}_{\lambda }^{\prime }(u),u\rangle =0\right\} .
\end{equation*}%
Clearly, $J_{\lambda }\left( u\right) $ is bounded below and coercive on $N$%
, since $p>2$. Following the standard argument, we conclude that there
exists a positive constant $\Lambda _{1}\geq \Lambda _{N}$ such that Eq. $(%
\ref{6-1})$ admits a positive ground state solution $\phi _{\lambda }\in
H^{2}(R^{N})$ for every $\lambda \geq \Lambda _{1}.$ Similar to the argument
of \cite[Theorem 22]{LWW}, we obtain that $\Pi _{\lambda }$ defined as (\ref%
{28}) is achieved and
\begin{equation*}
\Pi _{\lambda }=\frac{\left( \int_{\mathbb{R}^{N}}\left\vert \phi _{\lambda
}\right\vert ^{p}dx\right) ^{1/p}}{\left\Vert \phi _{\lambda }\right\Vert
_{\lambda }}>0\text{ for every }\lambda \geq \Lambda _{1}.
\end{equation*}%
Furthermore, similar to the proof of Theorem \ref{t1-3}, it follows that $%
\phi _{\lambda }\rightarrow \phi _{\infty }$ in $H^{2}(R^{N})$ and in $%
L^{p}(R^{N})$ as $\lambda \rightarrow \infty ,$ where $0\neq \phi _{\infty
}\in H_{0}^{1}(\Omega )\cap H^{2}(\Omega )$ is the weak solution of
biharmonic equations as follows%
\begin{equation*}
\Delta ^{2}u=|u|^{p-2}u\text{ in }\Omega .
\end{equation*}%
This implies that%
\begin{equation*}
\Pi _{\lambda }\rightarrow \Pi _{\infty }:=\frac{\left( \int_{\Omega
}\left\vert \phi _{\infty }\right\vert ^{p}dx\right) ^{1/p}}{\left(
\int_{\Omega }|\Delta \phi _{\infty }|^{2}dx\right) ^{1/2}}>0\text{ as }%
\lambda \nearrow \infty .
\end{equation*}%
Note that
\begin{equation*}
\Pi _{\lambda }=\sup_{u\in X_{\lambda }\backslash \{0\}}\frac{\left( \int_{%
\mathbb{R}^{N}}\left\vert u\right\vert ^{p}dx\right) ^{1/p}}{\left\Vert
u\right\Vert _{\lambda }}
\end{equation*}%
is decreasing on $\lambda .$ Hence, we have%
\begin{equation*}
\Pi _{\lambda }\searrow \Pi _{\infty }\text{ as }\lambda \nearrow \infty .
\end{equation*}

\section{Acknowledgments}

J. Sun was supported by the National Natural Science Foundation of China
(Grant No. 11671236), Shandong Provincial Natural Science Foundation (Grant
No. ZR2015JL002). T.F. Wu was supported in part by the Ministry of Science
and Technology (Grant No. 106-2115-M-390-001-MY2) and the National Center
for Theoretical Sciences, Taiwan.

\end{document}